\documentclass[reqno,12pt]{amsart}
\usepackage[utf8]{inputenc}
\usepackage[english]{babel}
\usepackage[symbol]{footmisc}
\usepackage{amssymb,amsmath,amsthm,amsfonts,xcolor,enumerate,comment}
\usepackage{cite}
\usepackage[mathcal]{euscript}
\usepackage{url}
\usepackage{geometry}
\usepackage{hyperref}

\theoremstyle{plain}
\newtheorem{theorem}{Theorem}[section]
\newtheorem{proposition}[theorem]{Proposition}
\newtheorem{lemma}[theorem]{Lemma}
\newtheorem{corollary}[theorem]{Corollary}

\theoremstyle{definition}
\newtheorem{definition}[theorem]{Definition}
\newtheorem{example}[theorem]{Example}

\theoremstyle{remark}
\newtheorem{remark}[theorem]{Remark}

\newcommand{\Hom}{\mathrm{Hom}}
\newcommand{\LL}{\mathfrak{L}}
\newcommand{\NN}{\mathfrak{N}}
\newcommand{\HH}{\mathfrak{H}}
\newcommand{\VV}{\mathfrak{V}}
\newcommand{\WW}{\mathfrak{W}}
\newcommand{\HV}{\mathfrak{HV}}
\newcommand{\sltwo}{\mathfrak{sl}_2}
\newcommand{\bF}{\mathbb{F}}
\newcommand{\bZ}{\mathbb{Z}}
\newcommand{\bN}{\mathbb{N}}
\newcommand{\bR}{\mathbb{R}}
\newcommand{\id}{\mathrm{id}}
\newcommand{\End}{\mathrm{End}}

\newtheorem*{theoremA}{Theorem A}
\newtheorem*{theoremB}{Theorem B}
\newtheorem*{theoremC}{Theorem C}
\newtheorem*{theoremD}{Theorem D}
\newtheorem*{theoremE}{Theorem E}
\newtheorem*{theoremF}{Theorem F}
\newtheorem*{theoremG}{Theorem G}

\newcommand{\repourl}{\url{https://aristotle.harmonic.fun/projects/f20edaf7-7c92-4351-a71b-c969f25d8231}}

\begin{document}

\title[$\Hom$-Lie structures on Lie algebras]
{$\Hom$-Lie algebra structures on Lie algebras}

\author{Jobir Adashev}
\author{Majidkhon Azizov}
\author{Vasily Voronin}

\begin{abstract}
We give a systematic description of $\Hom$-Lie structures on several classes of Lie algebras and
of their relations with transposed Poisson structures.  Our starting point is the observation
that the left-hand side of the defining identity is an alternating trilinear map; this makes
$\mathcal{H}(\LL)$, the set of $\Hom$-Lie structures of $\LL$, a subspace of $\End(\LL)$ and
reduces its computation to increasing triples of basis vectors.  From this we derive, in every
finite dimension $n$, the lower bound $\dim \mathcal{H}(\LL)\geq n \dim \mathcal{C}(\LL)$ and the
codimension bound $\binom{n}{3}\dim [\LL,[\LL,\LL]]$, together with a refinement of the latter,
and we determine $\mathcal{H}(\LL)$ for nilpotent, oscillator and solvable algebras, for
$\sltwo \oplus A_{n-3}$, for the Witt, Virasoro and Heisenberg--Virasoro algebras, for all
three-dimensional algebras, for $M_n(\bF)$ and $T_n(\bF)$, for the Virasoro-like algebra and its
$q$-analogue, and for the planar Galilean conformal algebra.  By a result of Filippov, every
$\frac12$-derivation, and more generally every $\delta$-derivation with $\delta\neq0,1$, of a Lie
algebra is a $\Hom$-Lie structure; we observe that, more generally, for $\delta\neq 0,1$ the
$\delta$-derivations of an algebra of any variety defined by identities of degree at most three
(associative, Novikov, left-symmetric, Leibniz, \ldots) are $\Hom$-structures of the
corresponding $\Hom$-variety.  We show that the
inclusion of the $\frac12$-derivations into the $\Hom$-Lie structures is far from being an
equality by computing the graded $\Hom$-Lie structures of Witt-, Virasoro- and
Schr\"odinger-type algebras and of their central extensions, for which we give an explicit
lifting obstruction.  Along the way we correct two multiplication tables from the literature that
violate the Jacobi identity.
\end{abstract}

\maketitle

\medskip


\medskip


\medskip

\noindent \textbf{Formal verification}: a number of the results of this paper are formalized
and machine-checked in Lean~4 with the Mathlib library.  In this version of the paper the proofs
are omitted; the formal proofs, together with the \LaTeX{} sources of the paper, are available at
\repourl.

\medskip

\setcounter{tocdepth}{1}
\tableofcontents

\section{Introduction}\label{sec:intro}

$\Hom$-Lie algebras were introduced by Hartwig, Larsson and Silvestrov in the study of
$q$-deformations of the Witt and Virasoro algebras \cite{HLS06}: a $\Hom$-Lie algebra is an
anticommutative algebra $\LL$ endowed with a linear map $\varphi$ such that the $\Hom$-Jacobi
identity
$$[\varphi(x),[y,z]]+[\varphi(y),[z,x]]+[\varphi(z),[x,y]]=0$$
holds for all $x,y,z \in \LL$.  If the underlying algebra is a Lie algebra, such a map $\varphi$
is called a \emph{$\Hom$-Lie structure} on $\LL$; the scalar multiples of the identity map are
$\Hom$-Lie structures on every Lie algebra and are called \emph{trivial}.  The problem of
describing all $\Hom$-Lie structures on a given Lie algebra has attracted a considerable
attention.  So, Jin and Li studied $\Hom$-Lie structures on semi-simple Lie algebras over an
algebraically closed field of characteristic zero \cite{jl08}, and Xie, Jin and Liu obtained the
complete description in the finite-dimensional simple case \cite[Theorem 3.3]{xjl15} (restated
as \cite[Theorem 3.1]{xl17}): all $\Hom$-Lie structures of a finite-dimensional simple Lie
algebra are trivial, with the single exception of $\sltwo$, whose $\Hom$-Lie structures form a
$6$-dimensional space.  $\Hom$-Lie structures were further described on simple graded Lie
algebras of finite growth \cite{xl17} --- a class which contains the Witt algebra, whose
$\Hom$-Lie structures are exactly the ``shifts with constant coefficients''
\cite[Theorem 3.2]{xl17} --- on affine Kac--Moody algebras \cite{MZ18}, where all $\Hom$-Lie
structures turn out to be spanned by the identity map and the maps with central image, on
finite-dimensional
simple Lie superalgebras \cite{CL13,yl15} and on simple Lie superalgebras of vector fields
\cite{sl17,ysl14}; for the Heisenberg algebras, the \emph{multiplicative} $\Hom$-Lie structures
were determined in \cite[Theorem 3.3]{Saa}.  The same principle --- twisting each defining
identity of degree three by a linear map applied to its ``outer'' variable --- produces the
$\Hom$-associative algebras of Makhlouf and Silvestrov \cite{MS08} and the $\Hom$-Novikov
algebras of Yau \cite{Yau11}; we come back to these varieties in Subsection \ref{subsec:delta}.

The present paper is devoted to a systematic study of several further classes of Lie algebras,
and to the interaction between $\Hom$-Lie structures and transposed Poisson structures.  The
notion of a transposed Poisson algebra was introduced by Bai, Bai, Guo and Wu \cite{BBGW} as a
dual notion of a Poisson algebra, obtained by exchanging the roles played by the two
multiplications in the Leibniz rule.  The systematic computation of transposed Poisson structures
on a given Lie algebra was initiated by Ferreira, Kaygorodov and Lopatkin \cite{FKL21}, who
related them to the $\frac{1}{2}$-derivations of the algebra.  Since then transposed Poisson
structures have been described for many classes of algebras, among them Witt type algebras
\cite{KL23,LSh26}, generalized Witt algebras and Block Lie algebras and superalgebras
\cite{KKh23,KKhGW,GS26,HGS25}, Galilean and solvable Lie algebras \cite{KLZ23}, solvable and
perfect Lie algebras \cite{KKhu24}, solvable Lie algebras with filiform nilradical \cite{AAE24},
quasi-filiform Lie algebras of maximum length \cite{ADSS}, the Schr\"odinger algebras
\cite{YTK,WC26}, the planar Galilean conformal algebra \cite{WZ} and the conformal Galilei
algebras \cite{CLY26}, the Lie algebras of upper triangular matrices and, more generally, the Lie
incidence algebras \cite{KKhUT,KKhInc}, Virasoro-, Witt- and Schr\"odinger--Virasoro-type
algebras \cite{KKhS24a,KKhS24b,Sh24}, the generalized Heisenberg--Virasoro algebra of rank two
\cite{WJSZ26}, the loop and the super Heisenberg--Virasoro algebras \cite{YT24b,YT24a}, the
Virasoro-like algebra \cite{LLJ2} and its $q$-analogue \cite{LLJ1}; a general approach through
$\frac{1}{2}$-(bi)derivations is developed in \cite{YH22}, a $\Hom$-version of the problem is
considered in \cite{LJ26}, and Poisson structures in the same spirit are described in
\cite{KKh21}.  The algebraic and geometric classification of the three-dimensional transposed
Poisson algebras is given in \cite{BFOK}, transposed $\delta$-Poisson and $\delta$-Poisson
algebras are studied in \cite{AKS26,AKSa26}, and the surveys \cite{K24,K25} give an account of
the state of the art.  The two notions are related by the observation that each multiplication
operator of a transposed Poisson structure satisfies the $\Hom$-Jacobi identity, so that a Lie
algebra without non-trivial $\Hom$-Lie structures has no non-trivial transposed Poisson
structures (Proposition \ref{prop:homtr}).

A second, general tool is provided by Theorem \ref{thm:half}: every $\frac{1}{2}$-derivation of
a Lie algebra is a $\Hom$-Lie structure.  This is a particular case of a result of Filippov
\cite{Fil98}, who proved that for $\delta\neq0,1$ every $\delta$-derivation $\varphi$ of a Lie
algebra satisfies $[\varphi(x),[y,z]]+[\varphi(y),[z,x]]+[\varphi(z),[x,y]]=0$, that is, is a
$\Hom$-Lie structure in the terminology of the present paper.  Recall that a linear map $\varphi$ of an algebra $A$ is a
\emph{$\delta$-derivation} if $\varphi(xy)=\delta\big(\varphi(x)y+x\varphi(y)\big)$; the
$\frac12$-derivations are the case $\delta=\frac12$.  The $\delta$-derivations were introduced by
Filippov, who studied them for Lie algebras, prime Lie algebras and prime alternative and
Mal'tsev algebras \cite{Fil98,Fil99,Fil00}; they were then described for classical Lie
superalgebras \cite{Kay09,Kay10}, for Lie algebras and superalgebras in general
\cite{Zus10}, for simple and semisimple Jordan algebras and superalgebras and for superalgebras
of Jordan brackets \cite{K07,Kay10,Kay12b,ZhK12}, for semisimple structurable algebras
\cite{KO14} and for $n$-ary algebras \cite{Kay12a}, and a $\delta$-analogue of the first
Whitehead lemma was obtained in \cite{ZZ21,ZZ25} (see also \cite{Kay11,Kay13} and the
survey \cite{K24}); more recently, $\delta$-derivations of Schr\"odinger and Virasoro-like
algebras were computed in \cite{WC26,CY25,LLJ2}, and $\delta$-derivations of commutative
associative algebras were used to construct $\delta$-Novikov and (transposed) $\delta$-Poisson
algebras \cite{Kay25,AKL25,AKS26}.  We show in Theorem \ref{thm:delta} that Filippov's theorem,
and hence the inclusion of Theorem \ref{thm:half}, is a special case of a general phenomenon:
\emph{for $\delta\neq 0,1$,
every $\delta$-derivation of an algebra of any variety defined by identities of degree at most
three is a $\Hom$-structure of the corresponding $\Hom$-variety}.  In particular, besides
Filippov's result that every $\delta$-derivation ($\delta\neq0,1$) of a Lie algebra is a
$\Hom$-Lie structure \cite{Fil98}, every
$\delta$-derivation of an associative algebra is a $\Hom$-associative structure in the sense of
\cite{MS08}, and every $\delta$-derivation of a Novikov algebra is a $\Hom$-Novikov structure in
the sense of \cite{Yau11}; the same holds for left-symmetric, Leibniz, Zinbiel, bicommutative and
$\delta$-Novikov algebras.  This makes the extensive literature on
$\frac{1}{2}$-derivations and $\delta$-derivations directly available for the present problem, and at the same time
raises the question of how far the inclusion $\Delta(\LL)\subseteq \mathcal{H}(\LL)$ is from
being an equality.  Sections \ref{sec:wab}--\ref{sec:nfg} answer this question for several
families of Witt-, Virasoro- and Schr\"odinger-type algebras whose $\frac{1}{2}$-derivations and
transposed Poisson structures were recently determined: the algebras $\mathcal{W}(a,b)$ and
the deformed generalized Heisenberg--Virasoro algebras \cite{KKhS25,KKhS24b}, the deformative
Schr\"odinger--Witt algebras $\mathcal{W}(a,b,s)$ and the extended Schr\"odinger--Witt algebra
\cite{KKhS24a,Sh24}, and the not-finitely graded Witt and Heisenberg--Witt algebras $W_n(G)$ and
$HW_n(G)$ \cite{KKhS24a,KKhS24b}.  In every one of them the space of $\Hom$-Lie structures turns
out to be infinite-dimensional --- also in the cases in which all $\frac{1}{2}$-derivations, and
hence all transposed Poisson structures, are trivial.

Sections \ref{sec:cext}--\ref{sec:nov} complete the picture by passing to the central
extensions, which are the algebras actually studied in \cite{KKhS24b,Sh24}: the deformed
generalized Heisenberg--Virasoro algebras $g(G,\lambda)$ (in particular the twisted
Heisenberg--Virasoro algebra), the not-finitely graded algebras $\widehat{W}(G)$,
$\widetilde{W}(G)$ and $\widetilde{HW}(G)$, the extended Schr\"odinger--Virasoro algebra
$\hat{\mathfrak{so}}$ and the original deformative Schr\"odinger--Virasoro algebras
$\widetilde{L_{\lambda,\mu}}^{\,i}$.  The computation is reduced, once and for all, to the
centreless case plus a single explicit obstruction $\Omega_{\psi}$ (Lemma \ref{lem:clift}),
which we then evaluate for each family.  Two of the multiplication tables printed in the
literature turn out not to satisfy the Jacobi identity; we exhibit explicit counterexamples and
prove that the corrected tables do define Lie algebras (Propositions \ref{prop:hwjac},
\ref{prop:hwhatlie}, \ref{prop:l4jac} and \ref{prop:l4lie}).

Our starting point, in Section \ref{sec:prelim}, is the elementary but useful observation that
the left-hand side
$$E_{\varphi}(x,y,z)=[\varphi(x),[y,z]]+[\varphi(y),[z,x]]+[\varphi(z),[x,y]]$$
of the defining identity is an \emph{alternating} trilinear function of its three arguments
(Lemma \ref{lem:alt}).  Consequently the identity only has to be tested on strictly increasing
triples of basis vectors (Corollary \ref{cor:triples}), and the $\Hom$-Lie structures of a Lie
algebra $\LL$ always form a linear subspace of $\End(\LL)$ containing the scalar maps and all
maps with central image (Lemma \ref{lem:subspace}).  This is the tool used in all the subsequent
sections.

The main results of the present paper are the following theorems.

\begin{theoremA}[see Proposition \ref{prop:2step}, Proposition \ref{prop:centralizer}, Theorem
\ref{thm:dimn}, Theorem \ref{thm:dim3} and Theorem \ref{thm:sl2}]
Let $\LL$ be a Lie algebra over a field $\bF$ of characteristic zero.
\begin{enumerate}
\item If $\LL$ is two-step nilpotent, in particular if $\LL$ is a Heisenberg algebra of arbitrary
dimension, then every linear map on $\LL$ is a $\Hom$-Lie structure.  (For Heisenberg algebras
this should be compared with the description of the \emph{multiplicative} $\Hom$-Lie structures
in \cite[Theorem 3.3]{Saa}; see Remark \ref{rem:heisknown}.)
\item If $\dim \LL=n<\infty$, then
$$\dim \mathcal{H}(\LL) \geq \max\Big\{\, n \cdot \dim \mathcal{C}(\LL), \ \
n^2-\binom{n}{3}\dim \LL^3 \,\Big\},$$
where $\LL^3=[\LL,[\LL,\LL]]$ and $\mathcal{C}(\LL)$ is the centralizer of $[\LL,\LL]$ in
$\LL$; the first bound holds because every linear map $\LL \to \mathcal{C}(\LL)$ is a
$\Hom$-Lie structure.  In particular, if $\dim \LL=3$, then the $\Hom$-Lie structures of $\LL$
form a subspace of $\End(\LL)$ of dimension at least $6$.
\item The bound in $(2)$ is attained for $n=3$ on $\sltwo=\langle h,e,f\rangle$: writing
$\varphi(h)=a_1h+a_2e+a_3f$, $\varphi(e)=b_1h+b_2e+b_3f$ and $\varphi(f)=c_1h+c_2e+c_3f$, the map
$\varphi$ is a $\Hom$-Lie structure if and only if $b_2=c_3$, $a_2=2c_1$ and $a_3=2b_1$.  Part
$(3)$ is not new: it is \cite[Theorem 3.3]{xjl15} (see also \cite[Theorem 3.1]{xl17}), which we
recover here by a direct computation; see Remark \ref{rem:sl2known}.
\end{enumerate}
\end{theoremA}

\begin{theoremB}[see Theorem \ref{thm:osc} and Theorem \ref{thm:codim1}]
Let $\LL_{\lambda}$ be the $(2n+2)$-dimensional oscillator Lie algebra.  Then the $\Hom$-Lie
structures of $\LL_{\lambda}$ are precisely the maps listed in Theorem \ref{thm:osc}.
Let $\LL$ be a solvable Lie algebra with abelian nilpotent radical $\NN$ of codimension $1$ and
$\dim \NN \geq 2$.  Then a linear map on $\LL$ is a $\Hom$-Lie structure if and only if it
preserves $\NN$.
\end{theoremB}

\begin{theoremC}[see Theorem \ref{thm:witt} and Theorem \ref{thm:vir}]
Let $\WW$ be the Witt algebra and $\VV$ its universal central extension, the Virasoro algebra.
\begin{enumerate}
\item A linear map $\varphi$ on $\WW$ is a $\Hom$-Lie structure if and only if
$\varphi(e_i)=\sum_{j \in \bZ}\alpha_j e_{i+j}$ for all $i \in \bZ$, where only finitely many
coefficients $\alpha_j$ are non-zero and they do not depend on $i$.  Part $(1)$ is not new: it
is \cite[Theorem 3.2]{xl17}, recovered here by a direct computation; see Remark
\ref{rem:wittknown}.
\item A linear map $\varphi$ on $\VV$ is a $\Hom$-Lie structure if and only if
$\varphi(x)=\alpha x+f(x)c$ for some scalar $\alpha$ and some linear form $f$ on $\VV$, i.e. if
and only if $\varphi$ is scalar modulo the centre.  Consequently, every transposed Poisson
structure on $\VV$ is trivial.  Part $(2)$ is the exact analogue, for the Virasoro algebra, of
the description of the $\Hom$-Lie structures of an affine Kac--Moody algebra given in
\cite[Theorem 3]{MZ18}; the Virasoro algebra is not an affine Kac--Moody algebra and is not
covered by loc. cit., and the result is obtained here independently (see Remark
\ref{rem:virkm}).  The last assertion is \cite[Corollary 28]{FKL21}, recovered here from the
description of
$\mathcal{H}(\VV)$.
\end{enumerate}
\end{theoremC}

\begin{theoremD}[see Theorem \ref{thm:half}, Remark \ref{rem:halfstrict}, Theorem
\ref{thm:delta} and Corollary \ref{cor:deltavar}]
Every $\frac12$-derivation of a Lie algebra $\LL$ is a $\Hom$-Lie structure, i.e.
$\Delta(\LL)\subseteq \mathcal{H}(\LL)$; more generally, for $\delta\neq0,1$ every
$\delta$-derivation of a Lie algebra is a $\Hom$-Lie structure (Filippov \cite{Fil98}).
The inclusion may be strict, and the quotient may be
arbitrarily large.  More generally, let $\delta \neq 0,1$ and let $A$ be an algebra of a variety
defined by multilinear identities of degree at most three.  Then every $\delta$-derivation of $A$
is a $\Hom$-structure of the corresponding $\Hom$-variety; in particular, $\delta$-derivations of
Lie, associative, left-symmetric, Novikov, Leibniz, Zinbiel, bicommutative and
$\delta'$-Novikov algebras (the Lie case being Filippov's theorem) are $\Hom$-Lie, $\Hom$-associative, $\Hom$-left-symmetric,
$\Hom$-Novikov, $\Hom$-Leibniz, $\Hom$-Zinbiel, $\Hom$-bicommutative and $\Hom$-$\delta'$-Novikov
structures respectively.
\end{theoremD}

\begin{theoremE}[see Theorem \ref{thm:wab}, Theorem \ref{thm:schro}, Theorem \ref{thm:eschro},
Theorem \ref{thm:wng} and Theorem \ref{thm:hwng}]
Let $\Gamma$ be a non-trivial additive subgroup of $\bF$.
\begin{enumerate}
\item A graded shift of degree $t$ of the algebra $W(\Gamma;a,b)$, i.e. a map
$L_{\mu}\mapsto \alpha L_{\mu+t}+\beta I_{\mu+t}$, $I_{\mu}\mapsto \gamma I_{\mu+t}$, is a
$\Hom$-Lie structure if and only if $\gamma=\alpha$ and $\alpha t(b^2-1)=0$.
\item A graded shift of degree $t$ of the deformative Schr\"odinger--Witt algebra
$\mathcal{S}(\Gamma;a,b)$, i.e. a map $L_{\mu}\mapsto \alpha L_{\mu+t}+\beta I_{\mu+t}$,
$I_{\mu}\mapsto \gamma I_{\mu+t}$, $Y_{\mu}\mapsto \delta Y_{\mu+t}$, is a $\Hom$-Lie structure
if and only if $\gamma=\delta=\alpha$ and $\alpha t(b+1)=0$.
\item Every map $L_{\mu}\mapsto \beta I_{\mu+t}$, $I_{\mu},Y_{\mu}\mapsto 0$ is a $\Hom$-Lie
structure of the algebras in $(1)$ and $(2)$, for all values of $a$ and $b$; likewise every map
$L_m \mapsto \beta M_{m+t}$, $M_n,N_n,Y_{n+\frac12}\mapsto 0$ is a $\Hom$-Lie structure of the
extended Schr\"odinger--Witt algebra, which has no non-trivial $\frac12$-derivation.
\item Every translation invariant map
$L_{\alpha,i}\mapsto \sum_{d,m}c_{d,m}L_{\alpha+d,i+m}$ is a $\Hom$-Lie structure of the
not-finitely graded Witt algebra $W_n(G)$.  On the Heisenberg--Witt algebra $HW_n(G)$ with
$n \neq 0$, by contrast, a graded shift acting by the same non-zero scalar on the $L$'s and on
the $H$'s is a $\Hom$-Lie structure only if the shift is trivial, while all the maps
$L_{\alpha,i}\mapsto \beta H_{\alpha+d,i+m}$, $H_{\alpha,i}\mapsto 0$ are $\Hom$-Lie structures.
\end{enumerate}
\end{theoremE}

\begin{theoremF}[see Theorems \ref{thm:dghv}, \ref{thm:hatw}, \ref{thm:tildew},
\ref{thm:hwhat}, \ref{thm:esv}, \ref{thm:dsv} and \ref{thm:l4}]
Let $\widehat{\LL}=\LL \oplus Z$ be a central extension of a Lie algebra $\LL$ by a $2$-cocycle
$\omega$.  A linear map $\varphi$ of $\widehat{\LL}$ with $\varphi(Z)\subseteq Z$ is a $\Hom$-Lie
structure if and only if its $\LL$-part $\psi$ is a $\Hom$-Lie structure of $\LL$ and
$\sum_{\mathrm{cyc}}\omega(\psi(x),[y,z])=0$.  Evaluating this obstruction on graded shifts one
obtains, for the central extensions listed above:
\begin{enumerate}
\item for $g(G,\lambda)$: the shift $L_a \mapsto \alpha L_{a+t}+\beta I_{a+t}$,
$I_a \mapsto \gamma I_{a+t}$ is a $\Hom$-Lie structure if and only if $\gamma=\alpha$ and
$\alpha t=0$, and moreover $\beta t=0$ for $\lambda=1$ and $\beta=0$ for $\lambda=0$ and for
$\lambda=-2$ of rank $\geq 2$; in particular, on the twisted Heisenberg--Virasoro algebra all
graded shifts which are $\Hom$-Lie structures are scalar modulo the centre;
\item for $\widehat{W}(G)$: the shift $L_{\alpha,i}\mapsto cL_{\alpha+d,i+m}$ is a $\Hom$-Lie
structure if and only if $m \geq 1$ or $cd=0$, which reproduces exactly the description of
$\Delta(\widehat{W}(G))$ of \cite{KKhS24b};
\item for $\widetilde{W}(G)$: only the shifts with $cd=cm=0$, i.e. the maps which are scalar
modulo the centre;
\item for the extended Schr\"odinger--Virasoro algebra $\hat{\mathfrak{so}}$ and for the
algebras $\widetilde{L_{\lambda,\mu}}^{\,i}$, $i \in \{1,2,3,5\}$: the maps
$L_m \mapsto \beta M_{m+t}$ survive for every $t$, so that these algebras carry
infinite-dimensional spaces of $\Hom$-Lie structures modulo the centre --- although all of them,
with the single exception of $\widetilde{L_{1,\mu}}^{\,1}$, have no non-trivial
$\frac12$-derivation at all; for $\widetilde{L_{1,\mu}}^{\,4}$ only the single value $t=-2\mu$
survives.
\end{enumerate}
\end{theoremF}

\begin{theoremG}[see Theorem \ref{thm:dim3exact}, Theorem \ref{thm:gl}, Theorem \ref{thm:tn},
Theorem \ref{thm:vlike}, Theorem \ref{thm:qvlike} and Theorem \ref{thm:pgca}]
The following algebras, whose $\frac12$-derivations and transposed Poisson structures have
recently been described in the literature, have the following spaces of $\Hom$-Lie structures.
\begin{enumerate}
\item (Dimension three.)  For every three-dimensional Lie algebra $\LL$ the space
$\mathcal{H}(\LL)$ has codimension exactly $\dim \LL^3$ in $\End(\LL)$, i.e.
$\dim \mathcal{H}(\LL)=9-\dim \LL^3$; over an algebraically closed field of characteristic zero
the possible values are $9$ (two-step nilpotent algebras), $8$, $7$ and $6$ (only $\sltwo$).
Compare with the classification of the three-dimensional transposed Poisson algebras
\cite{BFOK}.
\item (Matrices.)  For the Lie algebra $M_n(\bF)$ of all $n \times n$ matrices over an
algebraically closed field of characteristic zero,
$\mathcal{H}(M_n(\bF))=\bF\,\id \oplus \Hom(M_n(\bF),Z)$ for $n \geq 3$, of dimension $n^2+1$,
and $\dim \mathcal{H}(M_2(\bF))=10$.  For the Lie algebra $T_n(\bF)$ of upper triangular matrices
and $n \geq 4$,
$$\mathcal{H}(T_n(\bF))=\Hom\big(T_n(\bF),\langle \delta,e_{1n}\rangle\big)\oplus
\langle \id,A,B,A',B'\rangle, \qquad \dim \mathcal{H}(T_n(\bF))=n(n+1)+5,$$
where $A,B,A',B'$ are the four explicit maps of Definition \ref{def:tnmaps}.  In both cases the
space of $\frac12$-derivations, computed in \cite{KKhUT}, is much smaller: of dimension $2$,
resp. $n+2$.
\item (Virasoro-like algebras.)  The Virasoro-like algebra $\mathcal{V}$, with basis
$\{L_m\}_{m \in \bZ^2\setminus\{0\}}$ and $[L_m,L_n]=(m_1n_2-m_2n_1)L_{m+n}$, and its
$q$-analogue $\mathcal{V}_q$ with $q$ not a root of unity, admit only trivial $\Hom$-Lie
structures:
$\mathcal{H}(\mathcal{V})=\mathcal{H}(\mathcal{V}_q)=\bF\,\id$.  Consequently they admit no
non-trivial transposed Poisson structure, which recovers results of \cite{LLJ1,LLJ2}.  This is
the opposite behaviour to that of the Witt algebra, for which
$\Delta(\WW)=\mathcal{H}(\WW)$ is infinite-dimensional.
\item (Planar Galilean conformal algebra.)  A graded shift of degree $t$ of
$\mathcal{G}(\Gamma)$ is a $\Hom$-Lie structure if and only if it is of the form
$L_{\mu}\mapsto \alpha L_{\mu+t}+\beta I_{\mu+t}+\gamma J_{\mu+t}$,
$H_{\mu}\mapsto \alpha H_{\mu+t}$, $I_{\mu}\mapsto \alpha I_{\mu+t}$,
$J_{\mu}\mapsto \alpha J_{\mu+t}$ with $\alpha t=0$.  In particular
$\mathcal{H}(\mathcal{G}(\Gamma))$ is infinite-dimensional, although all the
$\frac12$-derivations and all the transposed Poisson structures of $\mathcal{G}(\bZ)$ are
trivial \cite{WZ}.
\end{enumerate}
\end{theoremG}

Thus the Witt algebra has a wealth of $\Hom$-Lie structures, whereas its central extension has
essentially none; the same phenomenon, with the centre replaced by an abelian ideal, is observed
for the Heisenberg--Virasoro algebra $W(0,0)$, on which we construct a large family of $\Hom$-Lie
structures (Theorem \ref{thm:hv}) and show that the shifts of the Witt algebra are not among them
(Remark \ref{rem:hv}).  On the transposed Poisson side, $\sltwo$ shows that the implication of
Proposition \ref{prop:homtr} cannot be reversed: all transposed Poisson structures on $\sltwo$
are trivial, although it admits a $6$-dimensional space of $\Hom$-Lie structures (Corollary
\ref{cor:sl2tp}); the triviality of the transposed Poisson structures on $\sltwo$ is a
particular case of \cite[Corollary 9]{FKL21}, while the $6$-dimensional space of $\Hom$-Lie
structures is what makes the example a counterexample to the converse implication.  Some
directions for further work are collected in Section \ref{sec:further}.

Throughout, $\bF$ denotes the ground field, which is assumed to be of characteristic zero; in
Section \ref{sec:osc} the oscillator algebras are considered, as usual, over $\bF=\bR$.  All
Lie algebras and all maps between them are over $\bF$, and $\End(\LL)$ denotes the space of
$\bF$-linear maps $\LL \to \LL$.  Bases of infinite-dimensional algebras are always understood in
the algebraic sense, so that every element is a finite linear combination of basis vectors.

\section{$\Hom$-Lie structures and transposed Poisson structures}\label{sec:prelim}

\subsection{The identity (HOM)}

\begin{definition}\label{def:hom}
Let $(\LL,[\cdot,\cdot])$ be a Lie algebra and let $\varphi \colon \LL \to \LL$ be a linear map.
Then $(\LL,[\cdot,\cdot],\varphi)$ is a \emph{$\Hom$-Lie structure} on $(\LL,[\cdot,\cdot])$ if
\begin{equation}\label{HOM}
[\varphi(x),[y,z]]+[\varphi(y),[z,x]]+[\varphi(z),[x,y]]=0
\qquad \text{for all } x,y,z \in \LL .
\end{equation}
A $\Hom$-Lie structure is \emph{trivial} if $\varphi$ is a scalar multiple of the identity map.
\end{definition}

By abuse of language we shall say that the map $\varphi$ itself is a $\Hom$-Lie structure on
$\LL$.  It is convenient to give a name to the left-hand side of (\ref{HOM}): for
$\varphi \in \End(\LL)$ we put
\begin{equation}\label{Ephi}
E_{\varphi}(x,y,z)=[\varphi(x),[y,z]]+[\varphi(y),[z,x]]+[\varphi(z),[x,y]],
\qquad x,y,z \in \LL,
\end{equation}
so that $\varphi$ is a $\Hom$-Lie structure if and only if $E_{\varphi}$ vanishes identically.
Note that $E_{\varphi}(x,y,z)$ is linear in $\varphi$ as well.

\begin{lemma}\label{lem:alt}
For every $\varphi \in \End(\LL)$ the map $E_{\varphi}$ is trilinear and alternating: it is
invariant under cyclic permutations of $(x,y,z)$, it changes sign under each transposition of two
of its arguments, and it vanishes whenever two of its arguments are equal.
\end{lemma}

\begin{corollary}\label{cor:triples}
Let $\{e_i\}_{i \in I}$ be a basis of $\LL$ indexed by a linearly ordered set $I$.  A linear map
$\varphi \in \End(\LL)$ is a $\Hom$-Lie structure on $\LL$ if and only if
$E_{\varphi}(e_i,e_j,e_k)=0$ for all $i<j<k$.
\end{corollary}

\begin{lemma}\label{lem:subspace}
The $\Hom$-Lie structures on a Lie algebra $\LL$ form a linear subspace $\mathcal{H}(\LL)$ of
$\End(\LL)$.  It contains all scalar maps $\alpha\,\id$ and all linear maps whose image lies in
the centre of $\LL$.
\end{lemma}

\subsection{Transposed Poisson structures}

Recall, following \cite{BBGW}, that a \emph{transposed Poisson structure} on a Lie algebra
$(\LL,[\cdot,\cdot])$ is a
commutative associative multiplication $\cdot$ on $\LL$ satisfying the transposed Leibniz rule
\begin{equation}\label{TP}
2\, z\cdot [x,y] = [z\cdot x, y] + [x, z\cdot y], \qquad x,y,z \in \LL;
\end{equation}
it is called \emph{trivial} if the multiplication is identically zero.  For $a \in \LL$ we write
$\mathcal{M}_a \colon x \mapsto a \cdot x$ for the corresponding multiplication operator.

\begin{lemma}\label{lem:mult}
Let $\cdot$ be a transposed Poisson structure on a Lie algebra $(\LL,[\cdot,\cdot])$.  Then
$\mathcal{M}_a$ is a $\Hom$-Lie structure on $(\LL,[\cdot,\cdot])$ for every $a \in \LL$.
\end{lemma}

By (\ref{TP}) the operators $\mathcal{M}_a$ are exactly the $\frac12$-derivations arising from
the multiplication; the description of the $\frac12$-derivations of a Lie algebra is the standard
first step in the computation of its transposed Poisson structures \cite{FKL21,K25}.

\begin{proposition}\label{prop:homtr}
Let $(\LL,[\cdot,\cdot])$ be a Lie algebra with $\dim \LL \geq 2$, i.e. $\LL$ is not spanned by a
single element.  If $\LL$ admits no non-trivial $\Hom$-Lie structure, then $\LL$ admits no
non-trivial transposed Poisson structure.
\end{proposition}

\begin{remark}\label{rem:homtr}
The hypothesis $\dim \LL \geq 2$ cannot be omitted.  On a one-dimensional (abelian) Lie algebra
$\LL = \langle e \rangle$ every linear map is a scalar multiple of the identity, so all
$\Hom$-Lie structures are trivial, while the field multiplication $e \cdot e = e$ is a
non-trivial transposed Poisson structure.
\end{remark}

The implication of Proposition \ref{prop:homtr} cannot be reversed either; a counterexample is
given in Corollary \ref{cor:sl2tp} below.

\subsection{$\frac12$-derivations}\label{subsec:half}

Lemma \ref{lem:mult} is a special case of a general fact, due to Filippov \cite{Fil98} (who
proved it for all $\delta$-derivations with $\delta\neq0,1$, see Corollary \ref{cor:deltavar}),
which will be used repeatedly below to compare our results with the literature on transposed
Poisson structures.

\begin{definition}\label{def:half}
A linear map $\varphi \in \End(\LL)$ is a \emph{$\frac12$-derivation} of the Lie algebra $\LL$ if
$$\varphi([x,y])=\tfrac12\big([\varphi(x),y]+[x,\varphi(y)]\big) \qquad \text{for all }
x,y \in \LL .$$
The space of $\frac12$-derivations of $\LL$ is denoted by $\Delta(\LL)$.
\end{definition}

\begin{theorem}[Filippov \cite{Fil98}]\label{thm:half}
Every $\frac12$-derivation of a Lie algebra is a $\Hom$-Lie structure; that is,
$$\Delta(\LL) \subseteq \mathcal{H}(\LL)$$
for every Lie algebra $\LL$.
\end{theorem}

Since the multiplication operators $\mathcal{M}_a$ of a transposed Poisson structure are exactly
$\frac12$-derivations, Lemma \ref{lem:mult} is the special case $\varphi=\mathcal{M}_a$ of Theorem
\ref{thm:half}.  More importantly, Theorem \ref{thm:half} turns every computation of
$\Delta(\LL)$ --- and a great many such computations are available, see
\cite{FKL21,KKh23,KKhGW,KKhu24,KL23,KLZ23,KKhUT,KKhInc,KKhS24a,KKhS24b,KKhS25,Sh24,AAE24,ADSS,
GS26,CLY26,LSh26,WJSZ26,WC26,YT24b,YH22,YTK,WZ,LLJ1,LLJ2,K24,K25} --- into a supply of
$\Hom$-Lie structures on $\LL$.

\begin{remark}\label{rem:halfstrict}
The inclusion of Theorem \ref{thm:half} is in general strict, and very far from being an
equality.  For $\sltwo$ one has $\Delta(\sltwo)=\langle \id \rangle$, whereas
$\mathcal{H}(\sltwo)$ is $6$-dimensional by Theorem \ref{thm:sl2}.  Sections
\ref{sec:wab}--\ref{sec:nfg} below exhibit whole families of algebras with
$\Delta(\LL)=\langle \id \rangle$ and $\mathcal{H}(\LL)$ infinite-dimensional.
\end{remark}

\subsection{$\delta$-derivations and $\Hom$-structures in other varieties}\label{subsec:delta}

Filippov \cite{Fil98} proved that for $\delta\neq0,1$ every $\delta$-derivation of a Lie algebra
is a $\Hom$-Lie structure; Theorem \ref{thm:half} is the case $\delta=\frac12$.  Filippov's theorem
is in turn a particular case of a general and very simple principle, valid for all varieties of
algebras defined by identities of degree at most three and for all $\delta \neq 0,1$.  In this subsection $A$ is an arbitrary (not necessarily associative) algebra
over a field $\bF$ with multiplication $(x,y)\mapsto xy$.

\begin{definition}\label{def:delta}
Let $\delta \in \bF$.  A linear map $\varphi \in \End(A)$ is a \emph{$\delta$-derivation} of $A$
if
$$\varphi(xy)=\delta\big(\varphi(x)y+x\varphi(y)\big) \qquad \text{for all } x,y \in A .$$
\end{definition}

For $\delta=1$ these are the derivations, for $\delta=\frac12$ the $\frac12$-derivations of
Definition \ref{def:half}, and for $\delta=-1$ the antiderivations; see
\cite{Fil98,Fil99,Fil00,Kay09,Kay10,Kay11,Kay12a,Kay12b,Kay13,KO14,K07,ZhK12,Zus10} for the
description of $\delta$-derivations of many classes of algebras and superalgebras.

\begin{definition}\label{def:homvar}
A \emph{multilinear polynomial of degree three} is an expression
$$f(x_1,x_2,x_3)=\sum_{\sigma \in S_3}\Big(a_{\sigma}\,(x_{\sigma(1)}x_{\sigma(2)})x_{\sigma(3)}
+b_{\sigma}\,x_{\sigma(1)}(x_{\sigma(2)}x_{\sigma(3)})\Big), \qquad a_{\sigma},b_{\sigma}\in \bF .$$
For $\varphi \in \End(A)$ its \emph{$\Hom$-version} is obtained by applying $\varphi$ to the
outer variable of each monomial:
$$f^{\varphi}(x_1,x_2,x_3)=\sum_{\sigma \in S_3}\Big(a_{\sigma}\,(x_{\sigma(1)}x_{\sigma(2)})
\varphi(x_{\sigma(3)})+b_{\sigma}\,\varphi(x_{\sigma(1)})(x_{\sigma(2)}x_{\sigma(3)})\Big).$$
Let $\mathcal{V}$ be a variety of algebras defined by a family of multilinear identities of
degrees two and three.  A linear map $\varphi$ on an algebra $A \in \mathcal{V}$ is a
\emph{$\Hom$-$\mathcal{V}$ structure} on $A$ if $f^{\varphi}(x,y,z)=0$ for all $x,y,z\in A$ and
every defining identity $f$ of degree three (the identities of degree two, such as commutativity
or anticommutativity, are not twisted).
\end{definition}

For the variety of Lie algebras, with $f=(x_1x_2)x_3+(x_2x_3)x_1+(x_3x_1)x_2$, one has
$f^{\varphi}(x,y,z)=-E_{\varphi}(x,y,z)$, so the $\Hom$-$\mathcal{V}$ structures are exactly the
$\Hom$-Lie structures.  In the same way Definition \ref{def:homvar} gives the $\Hom$-associative
algebras $(xy)\varphi(z)=\varphi(x)(yz)$ of \cite{MS08}, the $\Hom$-left-symmetric and
$\Hom$-Novikov algebras
$$(xy)\varphi(z)-\varphi(x)(yz)=(yx)\varphi(z)-\varphi(y)(xz), \qquad
(xy)\varphi(z)=(xz)\varphi(y)$$
of \cite{Yau11}, and the $\Hom$-Leibniz algebras
$(xy)\varphi(z)=(xz)\varphi(y)+\varphi(x)(yz)$.

\begin{lemma}\label{lem:deltamon}
Let $\varphi$ be a $\delta$-derivation of $A$.  Then for all $a,b,c \in A$
$$\varphi\big((ab)c\big)=\delta^{2}\big((\varphi(a)b)c+(a\varphi(b))c+(ab)\varphi(c)\big)
+\delta(1-\delta)\,(ab)\varphi(c),$$
$$\varphi\big(a(bc)\big)=\delta^{2}\big(\varphi(a)(bc)+a(\varphi(b)c)+a(b\varphi(c))\big)
+\delta(1-\delta)\,\varphi(a)(bc).$$
Consequently, for every multilinear polynomial $f$ of degree three,
$$\varphi\big(f(x,y,z)\big)=\delta^{2}\big(f(\varphi(x),y,z)+f(x,\varphi(y),z)+f(x,y,\varphi(z))\big)
+\delta(1-\delta)\,f^{\varphi}(x,y,z).$$
\end{lemma}

\begin{theorem}\label{thm:delta}
Let $\delta \in \bF\setminus\{0,1\}$ and let $f$ be a multilinear identity of degree three of the
algebra $A$.  Then every $\delta$-derivation $\varphi$ of $A$ satisfies
$f^{\varphi}(x,y,z)=0$ for all $x,y,z\in A$.  In particular, if $A$ belongs to a variety
$\mathcal{V}$ defined by multilinear identities of degrees two and three, then every
$\delta$-derivation of $A$ is a $\Hom$-$\mathcal{V}$ structure on $A$.
\end{theorem}

\begin{corollary}\label{cor:deltavar}
Let $\delta \neq 0,1$.
\begin{enumerate}
\item \textup{(Filippov \cite{Fil98})} Every $\delta$-derivation of a Lie algebra $\LL$ is a
$\Hom$-Lie structure on $\LL$.  For $\delta=\frac12$ this is Theorem \ref{thm:half}.
\item Every $\delta$-derivation of an associative algebra is a $\Hom$-associative structure in the
sense of \cite{MS08}.
\item Every $\delta$-derivation of a left-symmetric (respectively Novikov) algebra is a
$\Hom$-left-symmetric (respectively $\Hom$-Novikov) structure in the sense of \cite{Yau11}.
\item Every $\delta$-derivation of a Leibniz, Zinbiel or bicommutative algebra, and of a
$\delta'$-Novikov algebra in the sense of \cite{Kay25,AKL25} for any $\delta'$, is a
$\Hom$-structure of the corresponding variety, in the sense of Definition \ref{def:homvar}.
\end{enumerate}
\end{corollary}

\begin{remark}\label{rem:deltaexc}
The hypothesis $\delta\neq 0,1$ cannot be dropped.  For $\delta=1$: the derivation
$\varphi=\mathrm{ad}_e$ of $\sltwo$ is given by $\varphi(h)=-2e$, $\varphi(e)=0$, $\varphi(f)=h$;
in the notation of Theorem \ref{thm:sl2} it has $a_2=-2$ and $c_1=1$, so $a_2\neq 2c_1$ and
$\varphi$ is not a $\Hom$-Lie structure.  For $\delta=0$ the $0$-derivations are just the linear
maps vanishing on $A^2$, and Lemma \ref{lem:deltamon} gives no information.  Theorem
\ref{thm:delta} also says nothing about identities of degree four or higher, such as the Jordan or
the Mal'tsev identity.  Finally, Theorem \ref{thm:delta} extends verbatim to even
$\delta$-superderivations of superalgebras and the super-versions of the identities, since an
even map does not change the parities of the arguments; thus the even $\delta$-superderivations
($\delta\neq0,1$) computed, for instance, in \cite{Kay09,Kay10,Zus10,YT24a,HGS25} provide
$\Hom$-structures on the corresponding Lie superalgebras.
\end{remark}

\begin{remark}\label{rem:deltalit}
By Corollary \ref{cor:deltavar}, every known description of the $\delta$-derivations of a Lie
algebra with $\delta\neq 0,1$ --- for example those of \cite{Fil98,Fil99,Zus10,WC26,CY25,LLJ2}
--- yields $\Hom$-Lie structures, and $\delta$-derivations ($\delta\neq0,1$) of commutative
associative algebras,
which are the main tool of the constructions of \cite{Kay25,AKS26}, are
$\Hom$-associative structures.
\end{remark}

\section{Algebras of small dimension and two-step nilpotent algebras}\label{sec:small}

The criterion of Corollary \ref{cor:triples} is particularly efficient when the algebra has few
basis vectors or a small derived subalgebra.  We begin with two cases in which the identity
(\ref{HOM}) imposes no condition at all.

\begin{corollary}\label{cor:dim2}
If $\dim \LL \leq 2$, then every linear map $\varphi \colon \LL \to \LL$ is a $\Hom$-Lie
structure, i.e. $\mathcal{H}(\LL)=\End(\LL)$.
\end{corollary}

\begin{proposition}\label{prop:2step}
Let $\LL$ be two-step nilpotent, i.e. $[\LL,[\LL,\LL]]=0$.  Then every linear map
$\varphi \colon \LL \to \LL$ is a $\Hom$-Lie structure.  In particular this holds for every
abelian Lie algebra and for every Heisenberg algebra
$$\HH(V,\omega)=V \oplus \langle z \rangle, \qquad [v,w]=\omega(v,w)z, \qquad [v,z]=0
\qquad (v,w \in V),$$
associated with a skew-symmetric bilinear form $\omega$ on a vector space $V$ of arbitrary,
possibly infinite, dimension.
\end{proposition}

\begin{remark}\label{rem:heisknown}
For Heisenberg algebras Proposition \ref{prop:2step} should be compared with
\cite[Theorem 3.3]{Saa}, where the $\Hom$-Lie algebras whose underlying Lie algebra is the
$(2m+1)$-dimensional Heisenberg algebra $\HH_m=\langle x_1,\dots,x_m,y_1,\dots,y_m,z\rangle$,
$[x_k,y_k]=z$, are determined.  There the twisting map $\alpha$ is required to be
\emph{multiplicative}, i.e. to be in addition an endomorphism of $\HH_m$, and the answer is that
the matrix of $\alpha$ in the above basis has the block shape
$$P=\begin{pmatrix} X & T & 0 \\ Z & Y & 0 \\ L & M & \lambda \end{pmatrix},
\qquad \begin{pmatrix} X & T \\ Z & Y \end{pmatrix} \ \ \lambda\text{-symplectic} ,$$
that is, $\alpha(z)=\lambda z$ and the map induced by $\alpha$ on $V=\langle x_1,\dots,y_m\rangle$
multiplies the symplectic form by $\lambda$.  The two statements are consistent and
complementary.  Indeed, by Proposition \ref{prop:2step} the $\Hom$-Jacobi identity (\ref{HOM})
imposes no condition at all on a linear map of a two-step nilpotent algebra, so that in the
Heisenberg case the entire content of \cite[Theorem 3.3]{Saa} is the multiplicativity
requirement $\alpha([u,v])=[\alpha(u),\alpha(v)]$: for $u,v \in V$ this reads
$\omega(u,v)\alpha(z)=\omega(\alpha(u),\alpha(v))z$, which is precisely the condition that
$\alpha(z)=\lambda z$ and that the $V$-block be $\lambda$-symplectic.  Without the
multiplicativity assumption --- which is the setting of Definition \ref{def:hom} and of the
present paper --- the $\Hom$-Lie structures of $\HH_m$ form the whole space $\End(\HH_m)$, of
dimension $(2m+1)^2$; and Proposition \ref{prop:2step} covers in addition the Heisenberg
algebras $\HH(V,\omega)$ of infinite dimension and all two-step nilpotent algebras.
\end{remark}

\begin{remark}\label{rem:heis}
Heisenberg algebras also carry many non-trivial transposed Poisson structures.  Let $f$ be a
non-zero linear form on $V$ and put
$$(v+az)\cdot(w+bz)=f(v)f(w)\,z, \qquad v,w \in V, \ a,b \in \bF .$$
This multiplication is commutative, and it is associative because both
$((v+az)\cdot(w+bz))\cdot(u+cz)$ and $(v+az)\cdot((w+bz)\cdot(u+cz))$ are multiples of $f(z)=0$.
The transposed Leibniz rule (\ref{TP}) holds as well: its left-hand side vanishes because
$[x,y] \in \langle z \rangle$ and $f(z)=0$, and its right-hand side vanishes because the products
$z'\cdot x$ and $z'\cdot y$ are central.  This is consistent with Proposition \ref{prop:homtr}:
by Proposition \ref{prop:2step} the algebra $\HH(V,\omega)$ has plenty of non-trivial $\Hom$-Lie
structures.
\end{remark}

\subsection{Two lower bounds in arbitrary dimension}\label{subsec:dimn}

Corollary \ref{cor:dim2} and Proposition \ref{prop:2step} are the two extreme cases of a general
estimate, valid in every finite dimension, which we now state.  Throughout this subsection
$\LL$ is a Lie algebra of finite dimension $n$ over $\bF$, and we write
$$\LL^2=[\LL,\LL], \qquad \LL^3=[\LL,\LL^2], \qquad
\mathcal{C}(\LL)=\{x \in \LL \ : \ [x,\LL^2]=0\}$$
for the derived subalgebra, the third term of the lower central series, and the centralizer of
the derived subalgebra.  The subspace $\mathcal{C}(\LL)$ contains the centre of $\LL$, and
$\mathcal{C}(\LL)=\LL$ if and only if $\LL^3=0$, i.e. if and only if $\LL$ is two-step nilpotent.

The first bound produces $\Hom$-Lie structures, the second one bounds the number of conditions
that identity (\ref{HOM}) can impose.

\begin{proposition}\label{prop:centralizer}
Let $\LL$ be an arbitrary Lie algebra.  Then every linear map
$\varphi \colon \LL \to \mathcal{C}(\LL)$ is a $\Hom$-Lie structure, i.e.
$$\Hom(\LL,\mathcal{C}(\LL)) \subseteq \mathcal{H}(\LL).$$
If $\dim \LL = n < \infty$, then consequently
$$\dim \mathcal{H}(\LL) \geq n \cdot \dim \mathcal{C}(\LL).$$
\end{proposition}

Proposition \ref{prop:2step} is the special case $\mathcal{C}(\LL)=\LL$ of Proposition
\ref{prop:centralizer}, and the statement of Lemma \ref{lem:subspace} that all maps with central
image are $\Hom$-Lie structures is the special case
$\varphi(\LL)\subseteq Z(\LL)\subseteq \mathcal{C}(\LL)$.

\begin{theorem}\label{thm:dimn}
Let $\LL$ be a Lie algebra of dimension $n$ with basis $\{e_1,\ldots,e_n\}$.  For
$1 \leq i<j<k \leq n$ put
$$V_{ijk}=\big[\LL,\ \langle [e_i,e_j],\,[e_j,e_k],\,[e_k,e_i] \rangle\big] \subseteq \LL^3 .$$
Then the $\Hom$-Lie structures of $\LL$ form the kernel of the linear map
$$\Phi \colon \End(\LL) \longrightarrow \bigoplus_{1 \leq i<j<k \leq n} V_{ijk}, \qquad
\Phi(\varphi)=\big(E_{\varphi}(e_i,e_j,e_k)\big)_{1 \leq i<j<k \leq n},$$
and therefore
$$\dim \mathcal{H}(\LL) \ \geq \ n^2-\sum_{1 \leq i<j<k \leq n} \dim V_{ijk}
\ \geq \ n^2-\binom{n}{3}\dim \LL^3 .$$
\end{theorem}

\begin{corollary}\label{cor:dimn}
For every Lie algebra $\LL$ of finite dimension $n$,
$$\dim \mathcal{H}(\LL) \ \geq \
\max\Big\{\, n \cdot \dim \mathcal{C}(\LL),\ \ n^2-\binom{n}{3}\dim \LL^3 \,\Big\}.$$
In particular:
\begin{enumerate}
\item if $n \leq 2$, then $\binom{n}{3}=0$ and $\mathcal{H}(\LL)=\End(\LL)$ (Corollary
\ref{cor:dim2});
\item if $\LL$ is two-step nilpotent, then $\LL^3=0$ and $\mathcal{H}(\LL)=\End(\LL)$
(Proposition \ref{prop:2step});
\item if $n=3$, then $\dim \LL^3 \leq 3$ and $\binom{3}{3}=1$, so
$\dim \mathcal{H}(\LL)\geq 9-3=6$ (Theorem \ref{thm:dim3} below).
\end{enumerate}
\end{corollary}

\begin{remark}\label{rem:dimn}
The two bounds of Corollary \ref{cor:dimn} are of a complementary nature, and each of them is
vacuous in the range in which the other one is informative.

The bound $n^2-\binom{n}{3}\dim \LL^3$ is positive only for algebras which are close to being
two-step nilpotent: if $\dim \LL^3=d \geq 1$, it says something only when
$\binom{n}{3}d<n^2$, i.e. for $n \leq 8$ when $d=1$, for $n \leq 6$ when $d=2$, and for
$n \leq 5$ when $d=3$.  Its refined form $n^2-\sum_{i<j<k}\dim V_{ijk}$, however, keeps its
content in every dimension, because the triples $(i,j,k)$ with
$[e_i,e_j]=[e_j,e_k]=[e_k,e_i]=0$ --- which are the vast majority for a sparse multiplication
table --- contribute nothing at all; Example \ref{ex:fil4} shows that the refined bound can be
attained.

The bound $n \cdot \dim \mathcal{C}(\LL)$, on the other hand, is informative exactly when the
derived subalgebra has a large centralizer; it gives nothing for algebras with
$\mathcal{C}(\LL)=0$, such as $\sltwo$ and, more generally, any semi-simple Lie algebra.
Example \ref{ex:sl2ab} exhibits, in every dimension $n \geq 3$, an algebra for which the first
bound misses the true value by exactly $6$ while the second one is vacuous as soon as
$n \geq 6$.
\end{remark}

\begin{example}\label{ex:fil4}
Let $\NN_4=\langle e_1,e_2,e_3,e_4\rangle$ be the four-dimensional filiform Lie algebra, with
multiplication
$$[e_1,e_2]=e_3, \qquad [e_1,e_3]=e_4$$
and all the remaining brackets of basis vectors equal to zero.  Then
$$\mathcal{H}(\NN_4)=\{\varphi \in \End(\NN_4) \ : \ \varphi(e_3),\varphi(e_4) \in
\langle e_2,e_3,e_4\rangle\}, \qquad \dim \mathcal{H}(\NN_4)=14,$$
and this is exactly the value given by the refined bound of Theorem \ref{thm:dimn}.
\end{example}

\begin{example}\label{ex:sl2ab}
Let $n \geq 3$, let $m=n-3$ and let $\LL_n=\sltwo \oplus A_m$ be the direct sum of $\sltwo$ and
of an abelian Lie algebra $A_m=\langle a_1,\ldots,a_m\rangle$ of dimension $m$.  Then
$$\dim \mathcal{H}(\LL_n)=m^2+3m+6=n^2-3n+6=n \cdot \dim \mathcal{C}(\LL_n)+6 .$$
\end{example}

The family of Example \ref{ex:sl2ab} shows that the bound of Proposition \ref{prop:centralizer}
is asymptotically sharp --- the quotient of the two sides tends to $1$ as $n \to \infty$ --- but
never attained for $n \geq 3$, the defect being always equal to
$\dim \mathcal{H}(\sltwo)=6$.

\subsection{Dimension three}\label{subsec:dim3}

We now turn to the three-dimensional case, in which Corollary \ref{cor:triples} leaves exactly
one triple to be tested.

\begin{theorem}\label{thm:dim3}
Let $\LL$ be a three-dimensional Lie algebra with basis $\{e_1,e_2,e_3\}$.  Then the $\Hom$-Lie
structures of $\LL$ form the kernel of the linear map
$$\Phi \colon \End(\LL) \to \LL, \qquad \Phi(\varphi)=E_{\varphi}(e_1,e_2,e_3),$$
and therefore a subspace of $\End(\LL)$ of dimension at least $9-3=6$.
\end{theorem}

Theorem \ref{thm:dim3} is the case $n=3$ of Theorem \ref{thm:dimn}: there is a single
increasing triple, and its subspace $V_{123}$ is contained in the at most three-dimensional space
$\LL$.  The bound of Theorem \ref{thm:dim3} is attained, for instance by the simple Lie algebra
$\sltwo = \langle h,e,f \rangle$ with multiplication
\begin{equation}\label{sl2brackets}
[h,e]=2e, \qquad [h,f]=-2f, \qquad [e,f]=h .
\end{equation}

\begin{theorem}[Xie, Jin and Liu, {\cite[Theorem 3.3]{xjl15}}; see also
{\cite[Theorem 3.1]{xl17}}]\label{thm:sl2}
Write a linear map $\varphi$ on $\sltwo$ as
$$\varphi(h)=a_1h+a_2e+a_3f, \qquad \varphi(e)=b_1h+b_2e+b_3f, \qquad
\varphi(f)=c_1h+c_2e+c_3f .$$
Then $\varphi$ is a $\Hom$-Lie structure if and only if
$$b_2=c_3, \qquad a_2=2c_1, \qquad a_3=2b_1 .$$
Consequently the $\Hom$-Lie structures of $\sltwo$ form a $6$-dimensional space, parametrized
bijectively by $(a_1,b_1,b_2,b_3,c_1,c_2) \in \bF^6$ through
$$\varphi(h)=a_1h+2c_1e+2b_1f, \qquad \varphi(e)=b_1h+b_2e+b_3f, \qquad
\varphi(f)=c_1h+c_2e+b_2f .$$
\end{theorem}

\begin{remark}\label{rem:sl2known}
Theorem \ref{thm:sl2} is not a new result: the space $\mathcal{H}(\sltwo)$ was computed by Xie,
Jin and Liu \cite[Theorem 3.3]{xjl15}, and their answer is reproduced as
\cite[Theorem 3.1]{xl17}, where a $\Hom$-Lie structure is written as its matrix in the standard
basis $\{X,Y,H\}=\{e,f,h\}$ of $\sltwo$:
$$\mathcal{H}(\sltwo)=\langle E_{12},\ E_{21},\ E_{33},\ E_{11}+E_{22},\ 2E_{13}+E_{32},\
2E_{23}+E_{31}\rangle .$$
This is exactly the space described above: the six listed matrices correspond, in the
parametrization of Theorem \ref{thm:sl2}, to $c_2=1$, $b_3=1$, $a_1=1$, $b_2=c_3=1$,
$(a_2,c_1)=(2,1)$ and $(a_3,b_1)=(2,1)$ respectively, and conversely the relations $b_2=c_3$,
$a_2=2c_1$, $a_3=2b_1$ say precisely that the matrix of $\varphi$ lies in their span.  We keep
the statement here because it follows from Corollary \ref{cor:triples} by a short computation,
which is independent of the structure theory of semi-simple Lie algebras, and it
is valid over an arbitrary field of characteristic zero (algebraic closedness, assumed in
\cite{xjl15,xl17}, is not needed); it also shows that the general lower bound of Theorem
\ref{thm:dim3} is sharp.
\end{remark}

\begin{corollary}[first assertion: {\cite[Corollary 9]{FKL21}}]\label{cor:sl2tp}
Every transposed Poisson structure on $\sltwo$ is trivial, while $\sltwo$ admits non-trivial
$\Hom$-Lie structures.  In particular the implication of Proposition \ref{prop:homtr} cannot be
reversed.
\end{corollary}

\begin{remark}\label{rem:sl2tpfkl}
The first assertion of Corollary \ref{cor:sl2tp} is not a new result: it is the special case
$\LL=\sltwo$ of \cite[Corollary 9]{FKL21}, where Ferreira, Kaygorodov and Lopatkin prove that a
complex semi-simple finite-dimensional Lie algebra admits no non-trivial transposed Poisson
structure.  Their argument combines \cite[Theorem 8]{FKL21} --- a Lie algebra all of whose
$\frac12$-derivations are trivial carries only trivial transposed Poisson structures --- with
Filippov's description of the $\frac12$-derivations of simple finite-dimensional Lie algebras
(see \cite{FKL21,K07}).  The derivation of Corollary \ref{cor:sl2tp} from Theorem
\ref{thm:sl2} is self-contained, does not use the structure theory of semi-simple Lie algebras, and is valid over
an arbitrary field of characteristic zero, whereas \cite{FKL21} works over $\mathbb{C}$.  What is
relevant for the present paper is the combination of this known fact with the second assertion:
$\sltwo$ has a $6$-dimensional space of $\Hom$-Lie structures by Theorem \ref{thm:sl2}, so it
witnesses that the implication of Proposition \ref{prop:homtr} is strict.
\end{remark}

\begin{remark}\label{rem:sl2dim}
On $\sltwo$ the trivial $\Hom$-Lie structures $\varphi=\alpha\,\id$ form a $1$-dimensional
subspace of the $6$-dimensional space of all $\Hom$-Lie structures, so the bound of Theorem
\ref{thm:dim3} is sharp.  By contrast, on the three-dimensional Heisenberg algebra the space of
$\Hom$-Lie structures is the whole $9$-dimensional space $\End(\LL)$, by Proposition
\ref{prop:2step}.  Both extremes allowed by Theorem \ref{thm:dim3} therefore occur.
\end{remark}

\subsection{The exact dimension in dimension three}\label{subsec:dim3exact}

For a three-dimensional Lie algebra the bound of Theorem \ref{thm:dimn} is not merely a bound:
it is an equality, and the answer depends on nothing but the dimension of $\LL^3$.

\begin{theorem}\label{thm:dim3exact}
Let $\LL$ be a three-dimensional Lie algebra over an arbitrary field.  Then the space
$\mathcal{H}(\LL)$ of $\Hom$-Lie structures of $\LL$ has codimension exactly $\dim \LL^3$ in
$\End(\LL)$; that is,
$$\dim \mathcal{H}(\LL)=9-\dim \LL^3, \qquad \LL^3=[\LL,[\LL,\LL]] .$$
\end{theorem}

\begin{corollary}\label{cor:dim3list}
Let $\LL$ be a three-dimensional Lie algebra over an algebraically closed field $\bF$ of
characteristic zero.  Then $\dim \mathcal{H}(\LL)$ is given by the following table, in which the
algebras are listed in the standard classification and only the non-zero brackets of basis
vectors are displayed.
$$\begin{array}{llcc}
\text{algebra} & \text{multiplication} & \dim \LL^3 & \dim \mathcal{H}(\LL)\\[2pt]
\hline
A_3 \ \text{(abelian)} & - & 0 & 9\\
\NN_3 \ \text{(Heisenberg)} & [e_1,e_2]=e_3 & 0 & 9\\
\mathfrak{r}_2 \oplus \bF & [e_1,e_2]=e_2 & 1 & 8\\
\mathfrak{r}_3 & [e_1,e_2]=e_2,\ [e_1,e_3]=e_2+e_3 & 2 & 7\\
\mathfrak{r}_{3,\lambda} \ (\lambda \neq 0) & [e_1,e_2]=e_2,\ [e_1,e_3]=\lambda e_3 & 2 & 7\\
\sltwo & [h,e]=2e,\ [h,f]=-2f,\ [e,f]=h & 3 & 6
\end{array}$$
In particular every three-dimensional Lie algebra admits non-trivial $\Hom$-Lie structures, the
value $9$ is attained exactly on the two-step nilpotent algebras, and the minimum $6$ is attained
exactly on $\sltwo$.
\end{corollary}

\begin{remark}\label{rem:dim3tp}
Corollary \ref{cor:dim3list} is the $\Hom$-Lie counterpart of the classification of the
three-dimensional transposed Poisson algebras obtained by Beites, Fern\'andez Ouaridi and
Kaygorodov \cite{BFOK}.  The two pictures are quite different: the transposed Poisson structures
of a three-dimensional Lie algebra may well be all trivial --- this is the case for $\sltwo$ by
Corollary \ref{cor:sl2tp} --- while the space of $\Hom$-Lie structures is never smaller than
$6$-dimensional.  Theorem \ref{thm:dim3exact} also shows that in dimension three the
dimension of $\mathcal{H}(\LL)$ is determined by the single number $\dim \LL^3$; in particular
it is the same, namely $7$, for all the algebras $\mathfrak{r}_{3,\lambda}$ with
$\lambda \neq 0$ and for $\mathfrak{r}_3$.
\end{remark}

\section{$\Hom$-Lie structures of oscillator Lie algebras}\label{sec:osc}

\begin{definition}\label{def:osc}
The oscillator Lie algebra $\LL_{\lambda}$ is the $(2n+2)$-dimensional real Lie algebra with
canonical basis $\mathbb{B}=\{e_{-1},e_0,e_i,f_i\}_{i=1,\ldots,n}$ and multiplication
$$[e_{-1},e_i]=\lambda_i f_i, \qquad [e_{-1},f_i]=-\lambda_i e_i, \qquad [e_i,f_i]=e_0
\qquad (i=1,\ldots,n),$$
where $\lambda=(\lambda_1,\ldots,\lambda_n) \in \bR^n$ with
$0<\lambda_1 \leq \ldots \leq \lambda_n$; all the remaining brackets of basis elements vanish.
\end{definition}

\begin{theorem}\label{thm:osc}
Let $\LL_{\lambda}$ be the $(2n+2)$-dimensional oscillator Lie algebra, $n \geq 1$, and let
$\varphi$ be a linear map on $\LL_{\lambda}$.  Then $\varphi$ is a $\Hom$-Lie structure if and
only if it has the form
\begin{align*}
\varphi(e_i) &= \sum_{j=1}^{n}\alpha_{j,i}e_j+\sum_{j=1}^{n}\beta_{j,i}f_j+\delta_i e_0
  && (i=1,\ldots,n),\\
\varphi(f_i) &= \sum_{j=1}^{n}\xi_{j,i}e_j+\sum_{j=1}^{n}\eta_{j,i}f_j+\tau_i e_0
  && (i=1,\ldots,n),\\
\varphi(e_{-1}) &= \sum_{j=1}^{n}\alpha_{j,-1}e_j+\sum_{j=1}^{n}\beta_{j,-1}f_j
  +\gamma_{-1}e_{-1}+\delta_{-1}e_0, &&\\
\varphi(e_0) &= \delta_0 e_0, &&
\end{align*}
where the coefficients satisfy
$$\alpha_{j,i}=\frac{\lambda_i}{\lambda_j}\alpha_{i,j}, \qquad
\beta_{j,i}=\frac{\lambda_i}{\lambda_j}\xi_{i,j}, \qquad
\eta_{j,i}=\frac{\lambda_i}{\lambda_j}\eta_{i,j}, \qquad 1 \leq i,j \leq n,$$
and the coefficients $\delta_i$, $\tau_i$, $\delta_0$, $\gamma_{-1}$, $\alpha_{j,-1}$,
$\beta_{j,-1}$, $\delta_{-1}$ are arbitrary.
\end{theorem}

\begin{remark}\label{rem:osc}
The second relation of Theorem \ref{thm:osc} couples $\beta$ with $\xi$ and not, as the symmetry
of the first and third relations might suggest, with $\eta$: the triple $\{e_i,f_k,e_{-1}\}$
produces $\lambda_i\xi_{i,k}=\lambda_k\beta_{k,i}$, because $[\varphi(e_i),e_k]$ picks out the
$f$-component of $\varphi(e_i)$ while $[\varphi(f_k),f_i]$ picks out its $e$-component.  No
relation of the form $\beta_{j,i}=\frac{\lambda_i}{\lambda_j}\eta_{i,j}$ holds: for $n=1$ and
$\lambda_1=1$, the map $\varphi(e_1)=f_1$, $\varphi(f_1)=e_1$, $\varphi(e_{-1})=\varphi(e_0)=0$
satisfies the conditions of Theorem \ref{thm:osc} (here $\beta_{1,1}=\xi_{1,1}=1$ and
$\alpha_{1,1}=\eta_{1,1}=0$), so it is a $\Hom$-Lie structure, but $\beta_{1,1} \neq \eta_{1,1}$.
\end{remark}

\section{$\Hom$-Lie structures of solvable Lie algebras with abelian nilpotent radical of
codimension $1$}\label{sec:codim1}

\begin{definition}\label{def:codim1}
Let $\LL$ be the solvable Lie algebra with abelian nilpotent radical
$\NN=\langle e_i \rangle_{i \geq 2}$ of codimension $1$, spanned by the generators
$\{e_i\}_{i \in \bN}$ and subject to
$$[e_1,e_n]=e_n, \qquad n \geq 2 ,$$
all the remaining brackets of basis elements being zero.
\end{definition}

\begin{theorem}\label{thm:codim1}
Let $\LL$ be a solvable Lie algebra with abelian nilpotent radical $\NN$ of codimension $1$, and
assume $\dim \NN \geq 2$, which is the case in Definition \ref{def:codim1}.  Then a linear map
$\varphi$ on $\LL$ is a $\Hom$-Lie structure if and only if $\varphi$ preserves the radical,
$\varphi(\NN) \subseteq \NN$.  In the basis of Definition \ref{def:codim1} this reads
$$\varphi(e_1)=\sum_{i \in \bN}\alpha_{1,i}e_i, \qquad
\varphi(e_j)=\sum_{i \geq 2}\alpha_{j,i}e_i \quad (j>1),$$
the value $\varphi(e_1)$ being arbitrary.
\end{theorem}

\section{$\Hom$-Lie structures of the Witt algebra}\label{sec:witt}

\begin{definition}\label{def:witt}
Let $\WW$ be the Witt algebra, i.e. the Lie algebra with basis $\{e_i\}_{i \in \bZ}$ and
multiplication $[e_i,e_j]=(i-j)e_{i+j}$.
\end{definition}

\begin{theorem}[{\cite[Theorem 3.2]{xl17}}]\label{thm:witt}
Let $\WW$ be the Witt algebra and let $\varphi$ be a linear map on $\WW$.  Then $\varphi$ is a
$\Hom$-Lie structure if and only if
$$\varphi(e_i)=\sum_{j \in \bZ}\alpha_j e_{i+j}, \qquad i \in \bZ,$$
where only finitely many of the coefficients $\alpha_j$ are non-zero and, crucially, they do not
depend on $i$.  Equivalently, writing $\gamma=\alpha_0$,
$$\varphi(e_i)=\sum_{j \in \bZ \setminus \{0\}}\alpha_j e_{i+j}+\gamma e_i, \qquad i \in \bZ .$$
\end{theorem}

\begin{remark}\label{rem:wittknown}
Theorem \ref{thm:witt} is not a new result: it is \cite[Theorem 3.2]{xl17}.  There the algebra
$\WW$ is called the Virasoro algebra --- it is the centreless Virasoro algebra, with basis
$\{L_n\}_{n \in \bZ}$ and $[L_m,L_n]=(m-n)L_{m+n}$, i.e. our Witt algebra of Definition
\ref{def:witt} --- and the answer is stated in the form
$$\mathcal{H}(\WW)=\bigoplus_{i \in \bZ}\bF\sigma_i, \qquad \sigma_i(L_n)=L_{n+i} ,$$
which is exactly the description above: a finite linear combination $\sum_i \alpha_i \sigma_i$
is precisely a map $e_n \mapsto \sum_i \alpha_i e_{n+i}$ with finitely many non-zero
coefficients $\alpha_i$ not depending on $n$.  This is also the case of the ``two-sided Witt
algebra'' recorded in the introduction of \cite{MZ18}.  The proof in \cite{xl17} proceeds by
showing that $\mathcal{H}(\WW)$ is graded (the algebra possesses a finite commutant generating
set) and then treating the homogeneous components separately; our approach is a direct
computation with the identity (\ref{HOM}), and Theorem \ref{thm:witt} is valid over an arbitrary
field of characteristic zero (algebraic closedness, assumed in \cite{xl17}, is not used).
\end{remark}

\begin{remark}\label{rem:witt}
The independence of the coefficients $\alpha_j$ from $i$ is essential, and it is not implied by
the triples $\{e_0,e_i,e_j\}$ alone.  Indeed, those triples, evaluated on the specializations
$m=i-j$, $m=j-i$ and $m=0$, yield the relations
$$\alpha_{i,j}=\alpha_{0,j-i} \ \ (i \neq j), \qquad
\alpha_{i,i}=(i+1)\alpha_{0,0}-i\alpha_{1,1} \ \ (i \neq 1),$$
which would allow the diagonal coefficient to be $\gamma$ at $i=1$ and $(i+1)\beta-\gamma$ at
$i \neq 1$.  Such maps are not $\Hom$-Lie structures when $\beta \neq 0$.  For instance, let
$\beta=1$, $\gamma=0$ and $\alpha_j=0$ for all $j$, so that $\varphi(e_1)=0$ and
$\varphi(e_i)=(i+1)e_i$ for $i \neq 1$.  Then the triple $\{e_0,e_2,e_3\}$ gives
\begin{align*}
E_{\varphi}(e_0,e_2,e_3)
&=[\varphi(e_0),[e_2,e_3]]+[\varphi(e_2),[e_3,e_0]]+[\varphi(e_3),[e_0,e_2]]\\
&=[e_0,-e_5]+[3e_2,3e_3]+[4e_3,-2e_2]
=5e_5-9e_5-8e_5=-12\,e_5 \neq 0 .
\end{align*}
The triples $\{e_k,e_i,e_j\}$ with $k \neq 0$ thus impose genuinely stronger conditions, and by
Theorem \ref{thm:witt} they force $\beta=0$: no $\Hom$-Lie structure of the Witt algebra has a
diagonal coefficient depending on $i$.
\end{remark}

\section{$\Hom$-Lie structures of the Virasoro algebra}\label{sec:vir}

By Theorem \ref{thm:witt} the Witt algebra has a wealth of $\Hom$-Lie structures.  Its universal
central extension behaves in the opposite way: modulo the centre, the only $\Hom$-Lie structures
of the Virasoro algebra are the scalar ones.

\begin{definition}\label{def:vir}
Let $\VV$ be the Virasoro algebra, i.e. the Lie algebra with basis
$\{e_i\}_{i \in \bZ} \cup \{c\}$ and multiplication
$$[e_i,e_j]=(i-j)e_{i+j}+\delta_{i+j,0}\,\frac{i^3-i}{12}\,c, \qquad [e_i,c]=[c,e_i]=0 .$$
\end{definition}

\begin{theorem}\label{thm:vir}
Let $\VV$ be the Virasoro algebra and let $\varphi$ be a linear map on $\VV$.  Then $\varphi$ is a
$\Hom$-Lie structure if and only if there are a scalar $\alpha \in \bF$ and a linear form
$f \colon \VV \to \bF$ such that
$$\varphi(x)=\alpha x+f(x)\,c \qquad \text{for all } x \in \VV,$$
that is, if and only if $\varphi$ is a scalar multiple of the identity modulo the centre.
\end{theorem}

\begin{remark}\label{rem:vir}
The contrast with the Witt algebra is sharp.  The shift $\varphi(e_i)=e_{i+1}$ is a $\Hom$-Lie
structure of the Witt algebra by Theorem \ref{thm:witt} (take $\alpha_1=1$ and $\alpha_j=0$ for
$j \neq 1$), but the map $\varphi(e_i)=e_{i+1}$, $\varphi(c)=0$ is \emph{not} a $\Hom$-Lie
structure of the Virasoro algebra, by Theorem \ref{thm:vir}, since it is not scalar modulo the
centre.  It is exactly the central term $\delta_{i+j,0}\frac{i^3-i}{12}c$ of the Virasoro bracket
that rules out all the shifts $\alpha_j$ with $j \neq 0$.
\end{remark}

\begin{remark}\label{rem:virkm}
Theorem \ref{thm:vir} should be compared with \cite[Theorem 3]{MZ18}, which states that the
space of $\Hom$-Lie structures on an affine Kac--Moody algebra $\widehat{\mathfrak{g}}$ is
linearly spanned by the identity map and by the $\Hom$-Lie structures with central image, i.e.
by the maps $\widehat{\mathfrak{g}}\to \bF z$.  This is formally the same answer as the one
obtained here for the Virasoro algebra: $\varphi=\alpha\,\id+f(\cdot)c$.  The Virasoro algebra
is, however, not an affine Kac--Moody algebra --- it is the universal central extension of the
Witt algebra, not of a loop algebra --- and we have not found the Virasoro case stated in
\cite{MZ18}, whose Sections 4 and 5 treat the untwisted and the twisted affine Kac--Moody algebras.
Theorem \ref{thm:vir} is therefore the exact analogue of \cite[Theorem 3]{MZ18} for the
Virasoro algebra, proved here independently.  The mechanism behind the two results is the same,
and it is the one isolated in Lemma \ref{lem:clift} below: on the centreless algebra there is a
large supply of $\Hom$-Lie structures (Theorem \ref{thm:witt} for the Witt algebra), and passing
to a central extension all of them except the scalar ones are killed by the obstruction coming
from the defining $2$-cocycle.
\end{remark}

\begin{corollary}[{\cite[Corollary 28]{FKL21}}]\label{cor:virtp}
The Virasoro algebra admits no non-zero transposed Poisson structure.
\end{corollary}

Corollary \ref{cor:virtp} is likewise not new: it is \cite[Corollary 28]{FKL21}, obtained there
from the fact that the Virasoro algebra has no non-trivial $\frac12$-derivations
\cite[Theorem 27]{FKL21}, and the same conclusion for the Witt and the Virasoro algebras is
obtained in \cite{KL23} through a direct computation of their $\frac12$-derivations; here it is
deduced instead from the description of the $\Hom$-Lie structures given by Theorem
\ref{thm:vir}.

\section{$\Hom$-Lie structures of the Heisenberg--Virasoro algebra}\label{sec:hv}

Our last example shows that an abelian ideal can play the role played by the centre of the
Virasoro algebra in Remark \ref{rem:vir}.

\begin{definition}\label{def:hv}
Let $\HV=W(0,0)$ be the (centreless) Heisenberg--Virasoro algebra, i.e. the semidirect product of
the Witt algebra with an abelian ideal: it has basis
$\{L_i\}_{i \in \bZ} \cup \{I_i\}_{i \in \bZ}$ and multiplication
$$[L_i,L_j]=(i-j)L_{i+j}, \qquad [L_i,I_j]=-j\,I_{i+j}, \qquad [I_i,I_j]=0 .$$
\end{definition}

The element $I_0$ is central in $\HV$, and $\HV$ is the first member of the family
$W(a,b)=\WW \ltimes A_{a,b}$ mentioned in Section \ref{sec:further}.

\begin{theorem}\label{thm:hv}
For $t \in \bZ$ let $\psi_t$ be the linear map on $\HV$ with $\psi_t(L_i)=I_{i+t}$ and
$\psi_t(I_k)=0$.  Then $\psi_t$ is a $\Hom$-Lie structure of $\HV$, and hence so is every map
$$\varphi=\alpha\,\id+\sum_{t \in S}q_t\psi_t+f(\cdot)\,I_0$$
with $\alpha, q_t \in \bF$, $S \subseteq \bZ$ finite and $f$ an arbitrary linear form on $\HV$.
\end{theorem}

\begin{remark}\label{rem:hv}
As for the Virasoro algebra (Remark \ref{rem:vir}), the shifts of the Witt algebra do not survive:
the map $\varphi(L_i)=L_{i+1}$, $\varphi(I_k)=0$ is \emph{not} a $\Hom$-Lie structure of $\HV$.
Indeed, for the triple $(L_1,L_0,I_1)$,
\begin{align*}
E_{\varphi}(L_1,L_0,I_1)&=[L_2,[L_0,I_1]]+[L_1,[I_1,L_1]]+[0,[L_1,L_0]]\\
&=[L_2,-I_1]+[L_1,I_2]=I_3-2I_3=-I_3 \neq 0 .
\end{align*}
So only the ``vertical'' shifts $\psi_t$ of Theorem \ref{thm:hv}, which move the Witt part
\emph{into} the abelian ideal, remain.
\end{remark}

\section{$\Hom$-Lie structures of the algebras $W(a,b)$}\label{sec:wab}

The Heisenberg--Virasoro algebra of Section \ref{sec:hv} is the member $(a,b)=(0,0)$ of a
two-parameter family of semidirect products of the Witt algebra with an intermediate series
module.  It is no harder to treat these algebras in the generality of an arbitrary grading group,
which at the same time covers the centreless deformed generalized Heisenberg--Virasoro algebras.

\begin{definition}\label{def:wab}
Let $\Gamma$ be a non-trivial additive subgroup of $\bF$ and let $a,b \in \bF$.  The Lie algebra
$W(\Gamma;a,b)$ has basis $\{L_{\mu}\}_{\mu \in \Gamma} \cup \{I_{\mu}\}_{\mu \in \Gamma}$ and
multiplication
$$[L_{\mu},L_{\nu}]=(\mu-\nu)L_{\mu+\nu}, \qquad
[L_{\mu},I_{\nu}]=-(\nu+b\mu+a)I_{\mu+\nu}, \qquad [I_{\mu},I_{\nu}]=0 .$$
For $\Gamma=\bZ$ this is the algebra $W(a,b)=\WW \ltimes A_{a,b}$; thus
$W(\bZ;0,0)=\HV$ is the Heisenberg--Virasoro algebra of Definition \ref{def:hv}, and the
$L$-part is always a copy of the Witt algebra $\WW$ (of the higher rank Witt algebra
$W(\Gamma)$ in general).
\end{definition}

\begin{remark}\label{rem:wabdghv}
Two further families of the literature are members of this one.  The algebras
$\mathcal{W}(a,b)$ of \cite{KKhS25} are $W(\bZ;a,b)$.  The centreless deformed generalized
Heisenberg--Virasoro algebra $\check{g}(G,\lambda)$ of \cite{KKhS24b}, with
$[L_{\mu},L_{\nu}]=(\nu-\mu)L_{\mu+\nu}$ and $[L_{\mu},I_{\nu}]=(\nu-\lambda\mu)I_{\mu+\nu}$,
becomes $W(G;0,-\lambda)$ after the substitution $L_{\mu}\mapsto -L_{\mu}$.
\end{remark}

That $W(\Gamma;a,b)$ is a Lie algebra is checked on triples of basis vectors.  The only
non-obvious case is $(L_{\mu},L_{\nu},I_{\xi})$, for which, writing $u=\xi+a$,
\begin{align*}
&[L_{\mu},[L_{\nu},I_{\xi}]]+[L_{\nu},[I_{\xi},L_{\mu}]]+[I_{\xi},[L_{\mu},L_{\nu}]]\\
&\qquad =\Big((u+b\nu)(u+b\mu+\nu)-(u+b\mu)(u+b\nu+\mu)-(\mu-\nu)(u+b(\mu+\nu))\Big)I_{\mu+\nu+\xi}
\end{align*}
and the first two products differ by $(\nu-\mu)\big(u+b(\mu+\nu)\big)$, so the sum vanishes.

Since $W(\Gamma;a,b)$ is $\Gamma$-graded, with
$W(\Gamma;a,b)_{\mu}=\langle L_{\mu},I_{\mu}\rangle$, the natural class of maps to look at is
that of the homogeneous ones.

\begin{definition}\label{def:wabshift}
For $t \in \Gamma$ and $\alpha,\beta,\gamma \in \bF$ let $\varphi_{t;\alpha,\beta,\gamma}$ be the
\emph{graded shift} of degree $t$ determined by
$$\varphi_{t;\alpha,\beta,\gamma}(L_{\mu})=\alpha L_{\mu+t}+\beta I_{\mu+t}, \qquad
\varphi_{t;\alpha,\beta,\gamma}(I_{\mu})=\gamma I_{\mu+t} \qquad (\mu \in \Gamma).$$
\end{definition}

\begin{theorem}\label{thm:wab}
Let $\Gamma$ be a non-trivial additive subgroup of $\bF$ and let $a,b \in \bF$, $t \in \Gamma$,
$\alpha,\beta,\gamma \in \bF$.  Then $\varphi_{t;\alpha,\beta,\gamma}$ is a $\Hom$-Lie structure
of $W(\Gamma;a,b)$ if and only if
$$\gamma=\alpha \qquad \text{and} \qquad \alpha\,t\,(b^2-1)=0 .$$
In particular the coefficient $\beta$ is arbitrary.
\end{theorem}

\begin{corollary}\label{cor:wabpsi}
For all $t \in \Gamma$ and $\beta \in \bF$ the map
$$\psi_{t,\beta} \colon \quad L_{\mu} \mapsto \beta I_{\mu+t}, \qquad I_{\mu}\mapsto 0$$
is a $\Hom$-Lie structure of $W(\Gamma;a,b)$.  Consequently every map
$$\varphi=\alpha\,\id+\sum_{t \in S}\psi_{t,\beta_t}, \qquad S \subseteq \Gamma \ \text{finite},$$
is a $\Hom$-Lie structure, so $\mathcal{H}(W(\Gamma;a,b))$ is infinite-dimensional for every
$a,b$ whenever $\Gamma$ is infinite.
\end{corollary}

\begin{corollary}\label{cor:wabcases}
Let $t \in \Gamma$, $t \neq 0$.
\begin{enumerate}
\item If $b=1$ or $b=-1$, then $\varphi_{t;\alpha,\beta,\alpha}$ is a $\Hom$-Lie structure of
$W(\Gamma;a,b)$ for all $\alpha,\beta$.
\item If $b^2 \neq 1$, then the only graded shifts of non-zero degree that are $\Hom$-Lie
structures are the maps $\psi_{t,\beta}$ of Corollary \ref{cor:wabpsi}, whose $L$-part lands
inside the abelian ideal.
\end{enumerate}
\end{corollary}

\begin{remark}\label{rem:wabhv}
For $(a,b)=(0,0)$ and $\Gamma=\bZ$ we have $b^2-1=-1 \neq 0$, so Corollary \ref{cor:wabcases}(2)
recovers Theorem \ref{thm:hv} together with Remark \ref{rem:hv}: the maps $\psi_t=\psi_{t,1}$ are
$\Hom$-Lie structures of $\HV$, and the Witt-type shift $L_i \mapsto L_{i+1}$, $I_k \mapsto 0$
is not.
\end{remark}

\begin{remark}\label{rem:wabhalf}
By \cite{KKhS25} the algebra $\mathcal{W}(a,b)=W(\bZ;a,b)$ has non-trivial
$\frac12$-derivations only for $b=-1$.  Theorem \ref{thm:wab} shows that the $\Hom$-Lie picture
is completely different: the maps $\psi_{t,\beta}$ exist for \emph{all} values of $a$ and $b$, so
that $\Delta(W(\bZ;a,b))=\langle \id \rangle$ while $\mathcal{H}(W(\bZ;a,b))$ is
infinite-dimensional; moreover, non-trivial \emph{Witt-type} shifts survive not only for $b=-1$,
but also for $b=1$.  This is the announced strengthening of Remark \ref{rem:halfstrict}.
\end{remark}

\section{$\Hom$-Lie structures of the deformative Schr\"odinger--Witt algebras}\label{sec:schro}

We now add to the picture the third family of generators occurring in the Schr\"odinger-type
algebras.  Throughout this section we use the convention of \cite{KKhS24a,Sh24} for the Witt part,
i.e. $[L_{\mu},L_{\nu}]=(\nu-\mu)L_{\mu+\nu}$; it is turned into the convention of Section
\ref{sec:wab} by the substitution $L_{\mu}\mapsto -L_{\mu}$.

\begin{definition}\label{def:schro}
Let $\Gamma$ be a non-trivial additive subgroup of $\bF$ and let $a,b \in \bF$.  The
\emph{deformative Schr\"odinger--Witt algebra} $\mathcal{S}(\Gamma;a,b)$ has basis
$\{L_{\mu}\}\cup\{I_{\mu}\}\cup\{Y_{\mu}\}$, $\mu \in \Gamma$, and multiplication
\begin{align*}
[L_{\mu},L_{\nu}]&=(\nu-\mu)L_{\mu+\nu}, &
[L_{\mu},I_{\nu}]&=(\nu+b\mu+a)I_{\mu+\nu}, \\
[L_{\mu},Y_{\nu}]&=\Big(\nu+\tfrac{(b-1)\mu+a}{2}\Big)Y_{\mu+\nu}, &
[Y_{\mu},Y_{\nu}]&=(\nu-\mu)I_{\mu+\nu},
\end{align*}
all other brackets of basis vectors being zero.
\end{definition}

For $\Gamma=\bZ$ this is the algebra $\mathcal{W}(a,b,s)$ of \cite{KKhS24a} with $s=0$; since
$\mathcal{W}(a+1,b,0)\cong \mathcal{W}(a,b,\frac12)$, the half-integer case is covered as well.
The original deformative Schr\"odinger--Witt algebras $L_{\lambda,\mu}$ of \cite{Sh24} are the
members $\mathcal{W}(2\mu,-\lambda,\frac12)$ of the family.  The subalgebra spanned by the
$L$'s and the $I$'s is the algebra $W(\Gamma;a,b)$ of Section \ref{sec:wab} (up to the sign of
$L$), and the ideal $\mathcal{N}=\langle I_{\mu},Y_{\mu}\rangle$ is two-step nilpotent.

The Jacobi identity is verified on triples of basis vectors exactly as in Section \ref{sec:wab};
the two cases involving $Y$ reduce, after writing $c=\frac{b-1}{2}$ and $e=\frac{a}{2}$, to the
identities
$$(\xi+c\nu+e)(\nu+\xi+c\mu+e)-(\xi+c\mu+e)(\mu+\xi+c\nu+e)=(\nu-\mu)\big(\xi+e+c(\mu+\nu)\big)$$
for the triple $(L_{\mu},L_{\nu},Y_{\xi})$ and, for $(L_{\mu},Y_{\nu},Y_{\xi})$, to
$b\mu+a-\mu-2c\mu-2e=0$.

\begin{definition}\label{def:schroshift}
For $t \in \Gamma$ and $\alpha,\beta,\gamma,\delta \in \bF$ let
$\varphi_{t;\alpha,\beta,\gamma,\delta}$ be the graded shift of degree $t$ with
$$L_{\mu}\mapsto \alpha L_{\mu+t}+\beta I_{\mu+t}, \qquad
I_{\mu}\mapsto \gamma I_{\mu+t}, \qquad Y_{\mu}\mapsto \delta Y_{\mu+t}.$$
\end{definition}

\begin{theorem}\label{thm:schro}
Let $\Gamma$ be a non-trivial additive subgroup of $\bF$, $a,b \in \bF$, $t \in \Gamma$ and
$\alpha,\beta,\gamma,\delta \in \bF$.  Then $\varphi_{t;\alpha,\beta,\gamma,\delta}$ is a
$\Hom$-Lie structure of $\mathcal{S}(\Gamma;a,b)$ if and only if
$$\gamma=\delta=\alpha \qquad \text{and} \qquad \alpha\,t\,(b+1)=0 .$$
\end{theorem}

\begin{corollary}\label{cor:schropsi}
For all $t \in \Gamma$ and $\beta \in \bF$ the map
$L_{\mu}\mapsto \beta I_{\mu+t}$, $I_{\mu}\mapsto 0$, $Y_{\mu}\mapsto 0$ is a $\Hom$-Lie
structure of $\mathcal{S}(\Gamma;a,b)$, for all values of $a$ and $b$.  Hence
$\mathcal{H}(\mathcal{S}(\Gamma;a,b))$ is infinite-dimensional whenever $\Gamma$ is infinite.
\end{corollary}

\begin{remark}\label{rem:schrohalf}
By Theorem \ref{thm:schro} a graded shift with a non-zero Witt part and non-zero degree exists
exactly when $b=-1$ --- which is precisely the condition under which $\mathcal{W}(a,b,s)$ has
non-trivial $\frac12$-derivations, and hence non-trivial transposed Poisson structures, by
\cite{KKhS24a}.  In terms of the parameters of \cite{Sh24}, the algebra $L_{\lambda,\mu}$ is
exceptional exactly for $\lambda=1$.  Note, however, that even in the non-exceptional cases the
algebras carry the infinite family of $\Hom$-Lie structures of Corollary \ref{cor:schropsi}.
\end{remark}

\section{$\Hom$-Lie structures of the extended Schr\"odinger--Witt algebra}\label{sec:eschro}

The next algebra is the most rigid one from the point of view of transposed Poisson structures:
by \cite{Sh24} it has no non-trivial $\frac12$-derivations at all.  It nevertheless carries an
infinite-dimensional space of $\Hom$-Lie structures.

\begin{definition}[\cite{Sh24}]\label{def:eschro}
The \emph{extended Schr\"odinger--Witt algebra} $\widetilde{\mathfrak{so}}$ has basis
$\{L_n,M_n,N_n\}_{n \in \bZ}\cup\{Y_{n+\frac12}\}_{n \in \bZ}$ and multiplication
\begin{align*}
[L_m,L_n]&=(n-m)L_{m+n}, & [L_m,M_n]&=nM_{m+n}, & [L_m,N_n]&=nN_{m+n},\\
[N_m,M_n]&=2M_{m+n}, & [N_m,Y_{n+\frac12}]&=Y_{m+n+\frac12}, &
[Y_{m+\frac12},Y_{n+\frac12}]&=(m-n)M_{m+n+1},
\end{align*}
together with $[L_m,Y_{n+\frac12}]=\big(n+\frac{1-m}{2}\big)Y_{m+n+\frac12}$; all remaining
brackets of basis vectors vanish.
\end{definition}

\begin{theorem}\label{thm:eschro}
For every $t \in \bZ$ and every $\beta \in \bF$ the linear map $\chi_{t,\beta}$ determined by
$$\chi_{t,\beta}(L_m)=\beta M_{m+t}, \qquad
\chi_{t,\beta}(M_n)=\chi_{t,\beta}(N_n)=\chi_{t,\beta}(Y_{n+\frac12})=0$$
is a $\Hom$-Lie structure of $\widetilde{\mathfrak{so}}$.  Consequently
$$\alpha\,\id+\sum_{t \in S}\chi_{t,\beta_t}, \qquad S \subseteq \bZ \ \text{finite},$$
is a $\Hom$-Lie structure for all $\alpha,\beta_t \in \bF$, and
$\mathcal{H}(\widetilde{\mathfrak{so}})$ is infinite-dimensional.
\end{theorem}

\begin{remark}\label{rem:eschro}
By \cite{Sh24} every $\frac12$-derivation of $\widetilde{\mathfrak{so}}$ is trivial, so that
$\widetilde{\mathfrak{so}}$ carries no non-trivial transposed Poisson structure.  Theorem
\ref{thm:eschro} therefore provides a particularly striking illustration of
Remark \ref{rem:halfstrict}: here $\Delta(\widetilde{\mathfrak{so}})=\langle \id \rangle$ is
one-dimensional while $\mathcal{H}(\widetilde{\mathfrak{so}})$ is infinite-dimensional.  It also
shows that the converse of Proposition \ref{prop:homtr} fails for infinite-dimensional algebras
as badly as it does for $\sltwo$.
\end{remark}

\section{$\Hom$-Lie structures of not-finitely graded Witt and Heisenberg--Witt
algebras}\label{sec:nfg}

The last family we consider is the one of the not-finitely graded algebras $W_n(G)$ and
$HW_n(G)$, whose transposed Poisson structures were determined in \cite{KKhS24a}; their central
extensions $\widehat{W}(G)$, $\widetilde{W}(G)$ and $\widetilde{HW}(G)$ are treated in
\cite{KKhS24b}.  Here the grading group is $G \times \bZ$ and the bracket is no longer
homogeneous for the second component: it has a term of degree $0$ and a term of degree $n$.

\begin{definition}[\cite{KKhS24a}]\label{def:wng}
Let $G$ be a non-trivial additive subgroup of $\bF$ and let $n \in \bZ$.  The
\emph{not-finitely graded Witt algebra} $W_n(G)$ has basis
$\{L_{\alpha,i}\}_{\alpha \in G,\, i \in \bZ}$ and multiplication
$$[L_{\alpha,i},L_{\beta,j}]=(\beta-\alpha)L_{\alpha+\beta,i+j}+(j-i)L_{\alpha+\beta,i+j+n}.$$
The \emph{not-finitely graded Heisenberg--Witt algebra} $HW_n(G)$ has basis
$\{L_{\alpha,i},H_{\alpha,i}\}_{\alpha \in G,\, i \in \bZ}$, the same bracket on the $L$'s, and
$$[L_{\alpha,i},H_{\beta,j}]=\beta H_{\alpha+\beta,i+j}+jH_{\alpha+\beta,i+j+n}, \qquad
[H_{\alpha,i},H_{\beta,j}]=0 .$$
\end{definition}

\begin{theorem}\label{thm:wng}
For all $d \in G$, $m \in \bZ$ and $c \in \bF$ the graded shift
$$\sigma_{d,m,c} \colon \quad L_{\alpha,i}\longmapsto c\,L_{\alpha+d,i+m}$$
is a $\Hom$-Lie structure of $W_n(G)$.  Consequently every translation invariant map
$$\varphi(L_{\alpha,i})=\sum_{d \in G,\ m \in \bZ}c_{d,m}\,L_{\alpha+d,i+m}$$
with finitely many non-zero coefficients $c_{d,m}$, independent of $(\alpha,i)$, is a $\Hom$-Lie
structure of $W_n(G)$.
\end{theorem}

For $HW_n(G)$ the situation changes completely: the abelian ideal destroys the Witt-type shifts,
exactly as in Remark \ref{rem:hv}, while producing an infinite family of $\Hom$-Lie structures
with values in that ideal.

\begin{theorem}\label{thm:hwng}
Let $G$ be a non-trivial additive subgroup of $\bF$, let $n \in \bZ$, $d \in G$, $m \in \bZ$ and
$\alpha,\beta,\gamma \in \bF$, and let $\varphi$ be the graded shift of $HW_n(G)$ determined by
$$L_{\alpha_0,i}\mapsto \alpha L_{\alpha_0+d,i+m}+\beta H_{\alpha_0+d,i+m}, \qquad
H_{\alpha_0,i}\mapsto \gamma H_{\alpha_0+d,i+m}.$$
\begin{enumerate}
\item If $\gamma=\alpha$, $\gamma d=0$ and $\gamma m=0$, then $\varphi$ is a $\Hom$-Lie structure.
In particular $L_{\alpha_0,i}\mapsto \beta H_{\alpha_0+d,i+m}$, $H_{\alpha_0,i}\mapsto 0$ is a
$\Hom$-Lie structure for every shift $(d,m)$ and every $\beta$.
\item Conversely, if $n \neq 0$ and $\varphi$ is a $\Hom$-Lie structure with $\alpha=\gamma$,
then $\gamma d=0$ and $\gamma m=0$.
\end{enumerate}
\end{theorem}

\begin{remark}\label{rem:hwng}
For $n=0$ the condition for the graded shift $\varphi$ of Theorem \ref{thm:hwng} to be a
$\Hom$-Lie structure becomes $\alpha=\gamma$ and $\gamma(d+m)=0$; the surviving shifts are then those with $d \in G \cap \bZ$
and $m=-d$.  This matches exactly the description of $\Delta(HW_0(G))$ obtained in
\cite{KKhS24a}.
\end{remark}

\begin{remark}\label{rem:nfgcentral}
The algebras $\widehat{W}(G)$, $\widetilde{W}(G)$ and $\widetilde{HW}(G)$ of \cite{KKhS24b} are
central extensions of $W_1(G)$, $W_{-1}(G)$ and of a Heisenberg--Witt algebra.  As the passage
from the Witt algebra to the Virasoro algebra shows (Theorem \ref{thm:vir} and Remark
\ref{rem:vir}), a central term of the bracket can destroy all the graded shifts of Theorem
\ref{thm:wng}; this is settled in Section \ref{sec:nfgc}, where it turns out that the graded
shifts survive in $\widehat{W}(G)$ exactly when they move the second index upwards (Theorem
\ref{thm:hatw}), and that all the non-trivial ones die in $\widetilde{W}(G)$ (Theorem
\ref{thm:tildew}).
\end{remark}

\section{Central extensions: the obstruction to lifting a $\Hom$-Lie structure}\label{sec:cext}

All the algebras of Sections \ref{sec:wab}--\ref{sec:nfg} carry one-dimensional (or
finite-dimensional) central extensions which are the objects actually studied in
\cite{KKhS24b,Sh24}, and Remark \ref{rem:vir} shows that a central term of the bracket may
destroy all the graded shifts of the centreless algebra.  In this section we isolate the
mechanism behind this phenomenon; the computation of $\mathcal{H}$ for a central extension is
thereby reduced to the centreless case plus one explicit obstruction, which we then evaluate for
each of the families of \cite{KKhS24a,KKhS24b,Sh24}.

\begin{definition}\label{def:cext}
Let $(\LL,[\cdot,\cdot]_0)$ be a Lie algebra, let $Z$ be a vector space and let
$\omega \colon \LL \times \LL \to Z$ be an alternating bilinear map satisfying the cocycle
identity
\begin{equation}\label{cocycle}
\omega([x,y]_0,z)+\omega([y,z]_0,x)+\omega([z,x]_0,y)=0, \qquad x,y,z \in \LL .
\end{equation}
The \emph{central extension} $\widehat{\LL}=\LL \oplus Z$ is the Lie algebra with bracket
$$[x+u,\;y+v]=[x,y]_0+\omega(x,y), \qquad x,y \in \LL,\ u,v \in Z;$$
thus $Z$ is contained in the centre of $\widehat{\LL}$ and the projection
$\pi \colon \widehat{\LL}\to \LL$ along $Z$ is a homomorphism of Lie algebras.
\end{definition}

\begin{lemma}\label{lem:cproj}
Let $\varphi$ be a $\Hom$-Lie structure of $\widehat{\LL}$.  Then
$\overline{\varphi}:=\pi \circ \varphi|_{\LL}$ is a $\Hom$-Lie structure of $\LL$.
\end{lemma}

Conversely, a $\Hom$-Lie structure of $\LL$ lifts to $\widehat{\LL}$ precisely when one
obstruction vanishes.  For $\psi \in \End(\LL)$ we put
\begin{equation}\label{Omega}
\Omega_{\psi}(x,y,z)=\omega\big(\psi(x),[y,z]_0\big)+\omega\big(\psi(y),[z,x]_0\big)
+\omega\big(\psi(z),[x,y]_0\big) \in Z, \qquad x,y,z \in \LL .
\end{equation}

\begin{lemma}\label{lem:clift}
Every $\varphi \in \End(\widehat{\LL})$ can be written uniquely as
$$\varphi(x+u)=\psi(x)+f(x+u)+w(u) \qquad (x \in \LL,\ u \in Z)$$
with $\psi \in \End(\LL)$, $f \colon \widehat{\LL}\to Z$ linear and
$w \colon Z \to \LL$ linear.  Such a $\varphi$ is a $\Hom$-Lie structure of $\widehat{\LL}$ if
and only if
\begin{enumerate}
\item $\psi$ is a $\Hom$-Lie structure of $\LL$;
\item $\Omega_{\psi}=0$;
\item $[w(u),[y,z]]=0$ for all $u \in Z$ and all $y,z \in \widehat{\LL}$.
\end{enumerate}
In particular the ``central part'' $f$ is arbitrary, and if $w=0$ --- that is, if $\varphi(Z)
\subseteq Z$ --- then condition (3) is automatic.
\end{lemma}

\begin{corollary}\label{cor:cext0}
Let $\psi \in \End(\LL)$ and let $\mathrm{rad}(\omega)=\{x \in \LL : \omega(x,\LL)=0\}$.
\begin{enumerate}
\item If $\psi(\LL)\subseteq \mathrm{rad}(\omega)$, then $\Omega_{\psi}=0$.
\item $\Omega_{\alpha\,\id}=0$ for every scalar $\alpha$.
\item If $\omega$ is a coboundary, i.e. $\omega(x,y)=g([x,y]_0)$ for some linear
$g \colon \LL \to Z$, then $\Omega_{\psi}=0$ for every $\Hom$-Lie structure $\psi$ of $\LL$.
Consequently $\Omega_{\psi}$ only depends on the cohomology class of $\omega$, as long as $\psi$
is a $\Hom$-Lie structure of $\LL$.
\end{enumerate}
\end{corollary}

The following elementary remark will allow us to compute $\Omega$ on the not-finitely graded
algebras without any computation at all.

\begin{lemma}\label{lem:cpull}
Let $p \colon \LL \to A$ be a surjective homomorphism of Lie algebras, let $\omega_A$ be a
$Z$-valued $2$-cocycle on $A$ and let $\omega=p^{*}\omega_A$, i.e.
$\omega(x,y)=\omega_A(p(x),p(y))$; then $\omega$ is a $2$-cocycle on $\LL$.  If
$\psi \in \End(\LL)$ and $\psi_A \in \End(A)$ satisfy $p \circ \psi = \psi_A \circ p$, then
$$\Omega^{\omega}_{\psi}(x,y,z)=\Omega^{\omega_A}_{\psi_A}\big(p(x),p(y),p(z)\big),
\qquad x,y,z \in \LL .$$
In particular $\Omega^{\omega}_{\psi}=0$ whenever $\psi(\LL)\subseteq \ker p$.
\end{lemma}

All the obstructions computed below reduce to the following two polynomial identities.  We write
$$f(x)=x^3-x, \qquad g(p)=p(p-1)(p-2).$$

\begin{lemma}\label{lem:cubic}
Let $\bF$ be a field of characteristic zero.
\begin{enumerate}
\item For all $A,B,C \in \bF$,
$$(C-B)f(A)+(A-C)f(B)+(B-A)f(C)=(A-B)(B-C)(C-A)(A+B+C).$$
\item For all $u,v,w \in \bF$,
$$(w-v)g(u)+(u-w)g(v)+(v-u)g(w)=(u-v)(v-w)(w-u)(u+v+w-3).$$
\item If $\Gamma$ is a non-trivial additive subgroup of $\bF$ and $s \in \Gamma$, there are
$\alpha,\beta,\gamma \in \Gamma$ with $\alpha+\beta+\gamma=s$ and
$(\alpha-\beta)(\beta-\gamma)(\gamma-\alpha)\neq 0$.  The same holds for $\Gamma=\bZ$ and
$s \in \bZ$.
\item For all $m,n \in \bF$,
$$(n-m)f(m+n)+(m+2n)f(m)-(2m+n)f(n)=0 .$$
\end{enumerate}
\end{lemma}

The typical use of Lemma \ref{lem:cubic} is the following.  Let $\Gamma$ be a non-trivial
additive subgroup of $\bF$, let $\LL$ be $\Gamma$-graded with $\Gamma$-homogeneous elements
$L_{\mu}$ satisfying $[L_{\mu},L_{\nu}]_0=(\nu-\mu)L_{\mu+\nu}$, and let $\omega$ be a
\emph{Virasoro-type cocycle}, i.e. $\omega(L_{\mu},L_{\nu})=\delta_{\mu+\nu,0}\,f(\mu)\,C$ for a
central element $C$.  If $\psi(L_{\mu})=\alpha L_{\mu+t}+(\text{terms in }\mathrm{rad}\,\omega)$,
then for $\mu+\nu+\xi=-t$,
\begin{equation}\label{virobstr}
\Omega_{\psi}(L_{\mu},L_{\nu},L_{\xi})
=\alpha\sum_{\mathrm{cyc}}(\xi-\nu)f(\mu+t)\,C
=2\alpha t\,(\mu-\nu)(\nu-\xi)(\xi-\mu)\,C,
\end{equation}
by Lemma \ref{lem:cubic}(1) applied to $A=\mu+t$, $B=\nu+t$, $C=\xi+t$, whose sum is
$(\mu+\nu+\xi)+3t=2t$; all the other triples of basis vectors give $\Omega_{\psi}=0$ when
$\omega$ vanishes outside the $L$'s.  By Lemma \ref{lem:cubic}(3), $\Omega_{\psi}=0$ if and only
if $\alpha t=0$.  This is exactly the computation which, in Theorem \ref{thm:vir}, killed all the
shifts of the Witt algebra.

\section{$\Hom$-Lie structures of the deformed generalized Heisenberg--Virasoro
algebras}\label{sec:dghv}

We begin with the central extensions of the algebras $W(\Gamma;0,b)$ of Section \ref{sec:wab},
that is, with the deformed generalized Heisenberg--Virasoro algebras of \cite{KKhS24b}, whose
transposed Poisson structures were determined there.

\begin{definition}[\cite{KKhS24b}]\label{def:dghv}
Let $G$ be a non-trivial additive subgroup of $\bF$ and let $\lambda \in \bF$.  The
\emph{centreless deformed generalized Heisenberg--Virasoro algebra} $\check{g}(G,\lambda)$ has
basis $\{L_a,I_a\}_{a \in G}$ and multiplication
$$[L_a,L_b]=(b-a)L_{a+b}, \qquad [L_a,I_b]=(b-\lambda a)I_{a+b}, \qquad [I_a,I_b]=0 .$$
The \emph{deformed generalized Heisenberg--Virasoro algebra} $g(G,\lambda)$ is the central
extension $\check{g}(G,\lambda)\oplus Z$ of Definition \ref{def:cext} given by the cocycle
\begin{align*}
\omega(L_a,L_b)&=\delta_{a+b,0}\,\tfrac{a^3-a}{12}\,C_L, &
\omega(L_a,I_b)&=\delta_{a+b,0}\,h(a)\,C_{LI}, &
\omega(I_a,I_b)&=\delta_{a+b,0}\,a\,\delta_{\lambda,0}\,C_I,
\end{align*}
where
$$h(a)=\begin{cases}
a^2+a, & \lambda=0,\\
\frac{a^3-a}{12}, & \lambda=1,\\
\ell(a), & \lambda=-2,\\
0, & \text{otherwise},
\end{cases}$$
and, in the case $\lambda=-2$, $G$ is free of rank $\nu$ with a fixed $\bZ$-basis
$\varepsilon_1,\dots,\varepsilon_{\nu}$ and $\ell$ is one of the coordinate forms
$a \mapsto a_{(i)}$, $2 \leq i \leq \nu$ (for $\nu \geq 2$ there are $\nu-1$ independent central
elements $C^{(i)}_{LI}$; all the statements below are to be read for each of them separately).
The centre $Z$ is spanned by the central elements which actually occur.  For $G=\bZ$ and
$\lambda=0$ the algebra $g(\bZ,0)$ is the twisted Heisenberg--Virasoro algebra.
\end{definition}

By Remark \ref{rem:wabdghv} the substitution $L_a \mapsto -L_a$ identifies
$\check{g}(G,\lambda)$ with $W(G;0,-\lambda)$, so Theorem \ref{thm:wab} applies to it.  As in
Definition \ref{def:wabshift} we consider the graded shifts, now together with an arbitrary
central part.

\begin{definition}\label{def:dghvshift}
For $t \in G$, $\alpha,\beta,\gamma \in \bF$, a linear map $f \colon g(G,\lambda)\to Z$ and
elements $z_s \in Z$, let $\varphi$ be the linear map of $g(G,\lambda)$ given by
$$\varphi(L_a)=\alpha L_{a+t}+\beta I_{a+t}+f(L_a), \qquad
\varphi(I_a)=\gamma I_{a+t}+f(I_a), \qquad \varphi(C_s)=z_s$$
for the central generators $C_s$ of $Z$.  We call such a map a \emph{graded shift of degree $t$
of $g(G,\lambda)$}.
\end{definition}

\begin{lemma}\label{lem:thvobstr}
Let $A$ be the Lie algebra with basis $\{L_a,I_a\}_{a \in G}$ and
$[L_a,L_b]=(b-a)L_{a+b}$, $[L_a,I_b]=bI_{a+b}$, $[I_a,I_b]=0$ (that is,
$A=\check{g}(G,0)$), and let $\omega_A$ be the cocycle
$$\omega_A(L_a,L_b)=\delta_{a+b,0}\tfrac{a^3-a}{12}C_L, \quad
\omega_A(L_a,I_b)=\delta_{a+b,0}h(a)C_{LI}, \quad
\omega_A(I_a,I_b)=\delta_{a+b,0}\,a\,C_I$$
with $h(a)=a^2+\kappa a$ for some $\kappa \in \bF$.  Let
$\psi(L_a)=\alpha L_{a+t}+\beta I_{a+t}$ and $\psi(I_a)=\gamma I_{a+t}$ with $\gamma=\alpha$.
Then $\Omega_{\psi}=0$ if and only if $\beta=0$ and $\alpha t=0$.
\end{lemma}

\begin{theorem}\label{thm:dghv}
Let $\varphi$ be a graded shift of degree $t$ of $g(G,\lambda)$ as in Definition
\ref{def:dghvshift}.  Then $\varphi$ is a $\Hom$-Lie structure of $g(G,\lambda)$ if and only if
$$\gamma=\alpha, \qquad \alpha\,t=0,$$
together with
$$\beta\,t=0 \quad \text{if } \lambda=1, \qquad\qquad
\beta=0 \quad \text{if } \lambda=0, \text{ or if } \lambda=-2 \text{ and } \nu \geq 2 .$$
For all the remaining values of $\lambda$ the coefficient $\beta$ is arbitrary.
\end{theorem}

\begin{corollary}\label{cor:dghv}
Let $t \in G$.
\begin{enumerate}
\item If $\lambda \notin \{0,1,-2\}$, then $\psi_{t,\beta}\colon L_a \mapsto \beta I_{a+t}$,
$I_a \mapsto 0$, $C_s \mapsto 0$ is a $\Hom$-Lie structure of $g(G,\lambda)$ for every
$\beta \in \bF$; consequently $\mathcal{H}(g(G,\lambda))$ contains an infinite-dimensional space
of $\Hom$-Lie structures which are pairwise non-congruent modulo the maps with central image.
\item If $\lambda=0$, or if $\lambda=-2$ and $\nu \geq 2$, then every graded shift which is a
$\Hom$-Lie structure of $g(G,\lambda)$ is of the form $\alpha\,\id+(\text{a map into }Z)$ on
$\check{g}(G,\lambda)$.  In particular this holds for the twisted Heisenberg--Virasoro algebra
$g(\bZ,0)$.
\end{enumerate}
\end{corollary}

\begin{remark}\label{rem:dghv}
By \cite{KKhS24b} the algebras $g(G,\lambda)$ have non-trivial $\frac12$-derivations, and hence
non-trivial transposed Poisson structures, only for $\lambda=-1$.  Theorem \ref{thm:dghv} shows
that the $\Hom$-Lie picture is again different, but in a weaker way than for the centreless
algebras of Section \ref{sec:wab}: the vertical shifts $\psi_{t,\beta}$ survive for all
$\lambda \notin \{0,1,-2\}$, while for $\lambda \in \{0,-2\}$ --- in particular for the twisted
Heisenberg--Virasoro algebra --- the central terms $C_I$ and $C^{(i)}_{LI}$ destroy even those,
and only the trivial graded shifts remain.  Compare Theorem \ref{thm:hv} and Corollary
\ref{cor:wabpsi}, where $\psi_{t,\beta}$ was a $\Hom$-Lie structure for \emph{all} $a$ and $b$.
\end{remark}

\section{$\Hom$-Lie structures of the central extensions $\widehat{W}(G)$, $\widetilde{W}(G)$ and
$\widetilde{HW}(G)$}\label{sec:nfgc}

We now treat the three not-finitely graded central extensions of \cite{KKhS24b}; they are central
extensions of the algebras $W_n(G)$ and $HW_n(G)$ of Definition \ref{def:wng}.  Throughout,
$G$ is a non-trivial additive subgroup of $\bF$ and $\bZ_{\geq 0}$ denotes the set of
non-negative integers.  We write
$$W(G)_{+}=\langle L_{\alpha,i} \rangle_{\alpha \in G,\, i \in \bZ_{\geq 0}} \subseteq W_1(G),
\qquad
HW(G)_{+}=\langle L_{\alpha,i},H_{\alpha,i}\rangle_{\alpha \in G,\, i \in \bZ_{\geq 0}}
\subseteq HW_1(G);$$
these are subalgebras, because the second index of a bracket of two basis vectors of level
$\geq 0$ is again $\geq 0$.  The subspaces
$$\mathcal{I}_{+}=\langle L_{\alpha,i}\rangle_{i \geq 1}, \qquad
\mathcal{J}_{+}=\langle L_{\alpha,i},H_{\alpha,i}\rangle_{i \geq 1}$$
are ideals of $W(G)_{+}$ and of $HW(G)_{+}$ respectively, and the quotients are
\begin{equation}\label{quotients}
W(G)_{+}/\mathcal{I}_{+}\;\cong\;\WW(G), \qquad
HW(G)_{+}/\mathcal{J}_{+}\;\cong\;\check{g}(G,0),
\end{equation}
where $\WW(G)=\langle \bar{L}_{\alpha}\rangle$ with $[\bar L_{\alpha},\bar L_{\beta}]
=(\beta-\alpha)\bar L_{\alpha+\beta}$ is the higher rank Witt algebra and $\check{g}(G,0)$ is the
centreless twisted Heisenberg--Virasoro algebra of Lemma \ref{lem:thvobstr} (with
$\bar H_{\alpha}=\bar I_{\alpha}$).

\begin{definition}[\cite{KKhS24b}]\label{def:hatw}
The \emph{not-finitely graded Heisenberg--Virasoro algebra} $\widehat{W}(G)$ is the central
extension $W(G)_{+}\oplus \bF C$ given by the cocycle
$$\omega(L_{\alpha,i},L_{\beta,j})
=\delta_{\alpha+\beta,0}\,\delta_{i,0}\,\delta_{j,0}\,\tfrac{\alpha^3-\alpha}{12}\,C .$$
\end{definition}

That $\omega$ is a cocycle is immediate from (\ref{quotients}): $\omega$ is the pullback along
$p \colon W(G)_{+}\to \WW(G)$ of the Virasoro-type cocycle
$\omega_{\WW}(\bar L_{\alpha},\bar L_{\beta})=\delta_{\alpha+\beta,0}f(\alpha)C/12$ of
$\WW(G)$, which satisfies (\ref{cocycle}) by Lemma \ref{lem:cubic}(1): for
$\alpha+\beta+\gamma=0$,
$$\sum_{\mathrm{cyc}}\omega_{\WW}\big([\bar L_{\alpha},\bar L_{\beta}],\bar L_{\gamma}\big)
=\tfrac{1}{12}\sum_{\mathrm{cyc}}(\beta-\alpha)f(\alpha+\beta)\,C
=-\tfrac{1}{12}\sum_{\mathrm{cyc}}(\beta-\alpha)f(\gamma)\,C=0 .$$

\begin{theorem}\label{thm:hatw}
Let $d \in G$, $m \in \bZ_{\geq 0}$, $c \in \bF$, let $f \colon \widehat{W}(G)\to \bF C$ be a
linear map and $z \in \bF C$, and let $\varphi$ be the graded shift
$$\varphi(L_{\alpha,i})=c\,L_{\alpha+d,i+m}+f(L_{\alpha,i}), \qquad \varphi(C)=z .$$
Then $\varphi$ is a $\Hom$-Lie structure of $\widehat{W}(G)$ if and only if
$$m \geq 1 \qquad \text{or} \qquad c\,d=0 .$$
\end{theorem}

\begin{remark}\label{rem:hatwdelta}
By \cite{KKhS24b}, $\Delta(\widehat{W}(G))=\langle \id,\varphi\rangle$, where $\varphi$ runs
through the maps $L_{\alpha,i}\mapsto \sum_{d \in G}\sum_{k \geq 1}a^{d,k}L_{\alpha+d,i+k}$.
These are exactly the linear combinations of the graded shifts allowed by Theorem
\ref{thm:hatw}: here the graded part of $\mathcal{H}$ and $\Delta$ coincide, in sharp contrast
with the centreless algebras of Sections \ref{sec:wab}--\ref{sec:nfg}.  The reason is
structural: the shifts with $m \geq 1$ survive because they take values in the ideal
$\mathcal{I}_{+}$, on which the cocycle vanishes.
\end{remark}

\begin{definition}[\cite{KKhS24b}]\label{def:tildew}
The algebra $\widetilde{W}(G)$ is the central extension $W_{-1}(G)\oplus \bF C$ (Definition
\ref{def:wng} with $n=-1$, the second index running over $\bZ$) given by the cocycle
$$\omega(L_{\alpha,i},L_{\beta,j})=\delta_{\alpha+\beta,0}\Big(\delta_{i+j,-1}\alpha^3
+3i\,\delta_{i+j,0}\alpha^2+3i(i-1)\delta_{i+j,1}\alpha+i(i-1)(i-2)\delta_{i+j,2}\Big)C .$$
\end{definition}

\begin{theorem}\label{thm:tildew}
Let $d \in G$, $m \in \bZ$, $c \in \bF$, let $f \colon \widetilde{W}(G)\to \bF C$ be linear and
$z \in \bF C$, and let $\varphi$ be the graded shift
$$\varphi(L_{\alpha,i})=c\,L_{\alpha+d,i+m}+f(L_{\alpha,i}), \qquad \varphi(C)=z .$$
Then $\varphi$ is a $\Hom$-Lie structure of $\widetilde{W}(G)$ if and only if
$$c\,d=0 \qquad \text{and} \qquad c\,m=0,$$
that is, if and only if $\varphi$ is a scalar multiple of the identity modulo the centre.
\end{theorem}

\begin{remark}\label{rem:tildewdelta}
By \cite{KKhS24b}, $\Delta(\widetilde{W}(G))=\langle \id \rangle$ and $\widetilde{W}(G)$ has no
non-trivial transposed Poisson structure.  Theorem \ref{thm:tildew} shows that, for graded
shifts, the $\Hom$-Lie structures behave in the same way.  The contrast with Theorem
\ref{thm:hatw} is instructive: for $\widehat{W}(G)$ the cocycle is concentrated on the level
$i=0$, which is a quotient of the algebra, whereas the cocycle of $\widetilde{W}(G)$ involves
four levels, and the shifts in the $i$-direction can no longer avoid it.
\end{remark}

We turn to the last algebra of \cite{KKhS24b}, the one-dimensional central extension of a
Heisenberg--Witt algebra by three central elements.  Its multiplication table, as printed in
\cite[Section 5]{KKhS24b}, is
\begin{align}
[L_{\alpha,i},L_{\beta,j}]&=(\beta-\alpha)L_{\alpha+\beta,i+j}+(j-i)L_{\alpha+\beta,i+j-1}
+\delta_{\alpha,-\beta}\delta_{i,-j}\tfrac{\alpha^3-\alpha}{12}C_L, \label{hwt1}\\
[L_{\alpha,i},H_{\beta,j}]&=\beta H_{\alpha+\beta,i+j}+jH_{\alpha+\beta,i+j-1}
+\delta_{\alpha,-\beta}\delta_{i,-j}(\alpha^2-\alpha)C_{LH}, \label{hwt2}\\
[H_{\alpha,i},H_{\beta,j}]&=\delta_{\alpha,-\beta}\delta_{i,-j}\,\alpha\,C_H, \label{hwt3}
\end{align}
with $\alpha,\beta \in G$ and $i,j \in \bZ_{\geq 0}$.

\begin{proposition}\label{prop:hwjac}
The multiplication (\ref{hwt1})--(\ref{hwt3}) does not satisfy the Jacobi identity.  Explicitly,
$$\big[L_{-3,0},[L_{0,0},L_{3,1}]\big]+\big[L_{0,0},[L_{3,1},L_{-3,0}]\big]
+\big[L_{3,1},[L_{-3,0},L_{0,0}]\big]=-2\,C_L \neq 0 .$$
\end{proposition}

The failure is caused by the interaction of the cocycle, which is concentrated on the level
$i=0$, with the term $(j-i)L_{\alpha+\beta,i+j-1}$, which \emph{lowers} the level: the levels
$\geq 1$ do not form an ideal of $HW_{-1}(G)$, so a cocycle supported on the level $0$ cannot be
pulled back from a quotient.  If instead one uses the convention of $\widehat{W}(G)$ for the
Witt part, i.e. the algebra $HW_{1}(G)$ of Definition \ref{def:wng}, everything falls into place.

\begin{definition}\label{def:hwhat}
Let $\widetilde{HW}_1(G)=HW(G)_{+}\oplus \langle C_L,C_{LH},C_H\rangle$ be the central extension
of $HW(G)_{+}\subseteq HW_1(G)$ given by the cocycle
\begin{align*}
\omega(L_{\alpha,i},L_{\beta,j})&=\delta_{\alpha+\beta,0}\delta_{i,0}\delta_{j,0}
\tfrac{\alpha^3-\alpha}{12}C_L, &
\omega(L_{\alpha,i},H_{\beta,j})&=\delta_{\alpha+\beta,0}\delta_{i,0}\delta_{j,0}
(\alpha^2-\alpha)C_{LH}, \\
\omega(H_{\alpha,i},H_{\beta,j})&=\delta_{\alpha+\beta,0}\delta_{i,0}\delta_{j,0}\,\alpha\,C_H. &&
\end{align*}
\end{definition}

\begin{proposition}\label{prop:hwhatlie}
$\widetilde{HW}_1(G)$ is a Lie algebra, i.e. the above $\omega$ satisfies (\ref{cocycle}).
\end{proposition}

\begin{theorem}\label{thm:hwhat}
Let $d \in G$, $m \in \bZ_{\geq 0}$, $a,b,g \in \bF$, let $f$ be a linear map of
$\widetilde{HW}_1(G)$ into its centre $Z=\langle C_L,C_{LH},C_H\rangle$ and let
$z_L,z_{LH},z_H \in Z$.  Let $\varphi$ be the graded shift
$$\varphi(L_{\alpha,i})=a\,L_{\alpha+d,i+m}+b\,H_{\alpha+d,i+m}+f(L_{\alpha,i}), \qquad
\varphi(H_{\alpha,i})=g\,H_{\alpha+d,i+m}+f(H_{\alpha,i}),$$
$\varphi(C_L)=z_L$, $\varphi(C_{LH})=z_{LH}$, $\varphi(C_H)=z_H$.  Then $\varphi$ is a
$\Hom$-Lie structure of $\widetilde{HW}_1(G)$ if and only if
$$a=g, \qquad a\,d=0, \qquad a\,m=0, \qquad \text{and} \qquad b=0 \ \text{ if } m=0 .$$
\end{theorem}

\begin{remark}\label{rem:hwhat}
Thus $\widetilde{HW}_1(G)$ behaves like $\widehat{W}(G)$ in the $i$-direction and like the
twisted Heisenberg--Virasoro algebra in the $\alpha$-direction: the only graded shifts which are
$\Hom$-Lie structures and are non-trivial modulo the centre are the maps
$L_{\alpha,i}\mapsto b\,H_{\alpha+d,i+m}$, $H_{\alpha,i}\mapsto 0$ with $m \geq 1$, which take
their values in the ideal $\mathcal{J}_{+}$ on which the cocycle vanishes.  In particular
$\mathcal{H}(\widetilde{HW}_1(G))$ is infinite-dimensional modulo the maps with central image.
\end{remark}

\section{$\Hom$-Lie structures of the extended Schr\"odinger--Virasoro algebra}\label{sec:esv}

\begin{definition}[\cite{Sh24}]\label{def:esv}
The \emph{extended Schr\"odinger--Virasoro algebra} $\hat{\mathfrak{so}}$ is the central
extension of the extended Schr\"odinger--Witt algebra $\widetilde{\mathfrak{so}}$ of Definition
\ref{def:eschro} by $Z=\langle C_L,C_{LN},C_N\rangle$ given by the cocycle
$$\omega(L_m,L_n)=\delta_{m+n,0}\tfrac{m^3-m}{12}C_L, \quad
\omega(L_m,N_n)=\delta_{m+n,0}(m^2-m)C_{LN}, \quad
\omega(N_m,N_n)=\delta_{m+n,0}\,n\,C_N,$$
all the other values of $\omega$ on pairs of basis vectors being zero.
\end{definition}

\begin{theorem}\label{thm:esv}
For every $t \in \bZ$ and every $\beta \in \bF$ let $\hat\chi_{t,\beta}$ be the linear map of
$\hat{\mathfrak{so}}$ with
$$\hat\chi_{t,\beta}(L_m)=\beta M_{m+t}, \qquad
\hat\chi_{t,\beta}(M_n)=\hat\chi_{t,\beta}(N_n)=\hat\chi_{t,\beta}(Y_{n+\frac12})=0, \qquad
\hat\chi_{t,\beta}(Z)=0 .$$
Then $\hat\chi_{t,\beta}$ is a $\Hom$-Lie structure of $\hat{\mathfrak{so}}$, and so is every map
$$\varphi=\alpha\,\id+\sum_{t \in S}\hat\chi_{t,\beta_t}+f(\cdot), \qquad
S \subseteq \bZ \ \text{finite},$$
where $\alpha,\beta_t \in \bF$ and $f$ is an arbitrary linear map of $\hat{\mathfrak{so}}$ into
its centre.  In particular $\mathcal{H}(\hat{\mathfrak{so}})$ contains an infinite-dimensional
space of $\Hom$-Lie structures which are pairwise non-congruent modulo the maps with central
image.
\end{theorem}

\begin{remark}\label{rem:esv}
By \cite{Sh24} the algebra $\hat{\mathfrak{so}}$ has no non-trivial $\frac12$-derivation, and
hence no non-trivial transposed Poisson structure.  Together with Theorem \ref{thm:eschro} and
Remark \ref{rem:eschro} this gives a second, and even more rigid, example of the phenomenon of
Remark \ref{rem:halfstrict}: for both $\widetilde{\mathfrak{so}}$ and its central extension
$\hat{\mathfrak{so}}$ the space of $\frac12$-derivations is $\langle \id\rangle$, while the space
of $\Hom$-Lie structures is infinite-dimensional even modulo the centre.  In contrast with the
Virasoro algebra (Remark \ref{rem:vir}) and with the twisted Heisenberg--Virasoro algebra
(Corollary \ref{cor:dghv}(2)), passing to the central extension changes nothing here, because
the three cocycles of Definition \ref{def:esv} do not involve the abelian ideal $\mathcal{M}$ in
which the maps $\chi_{t,\beta}$ take their values.
\end{remark}

\section{$\Hom$-Lie structures of the original deformative Schr\"odinger--Virasoro
algebras}\label{sec:dsv}

The last family of \cite{Sh24} consists of the central extensions
$\widetilde{L_{\lambda,\mu}}^{\,i}$, $1 \leq i \leq 5$, of the original deformative
Schr\"odinger--Witt algebras $L_{\lambda,\mu}$, which are the members
$\mathcal{S}(\bZ;2\mu,-\lambda)$ of the family of Definition \ref{def:schro} (with the $Y$'s
indexed by $\frac12+\bZ$; see the discussion following that definition).  We recall the
multiplication tables, in the notation of Definition \ref{def:schro}: $L_{\lambda,\mu}$ has
basis $\{L_m,M_m,Y_{m+\frac12}\}_{m \in \bZ}$ and
\begin{align*}
[L_m,L_n]&=(n-m)L_{m+n}, & [L_m,M_n]&=(n-\lambda m+2\mu)M_{m+n},\\
[L_m,Y_{n+\frac12}]&=\Big(n+\tfrac12-\tfrac{\lambda+1}{2}m+\mu\Big)Y_{m+n+\frac12}, &
[Y_{m+\frac12},Y_{n+\frac12}]&=(n-m)M_{m+n+1} .
\end{align*}

\begin{definition}[\cite{Sh24}]\label{def:dsv}
The algebras $\widetilde{L_{\lambda,\mu}}^{\,i}$ are the central extensions of $L_{\lambda,\mu}$
determined by the following cocycles, where $f(x)=x^3-x$:
\begin{enumerate}
\item[$(1)$] $\omega(L_m,L_n)=\delta_{m+n,0}\frac{f(m)}{12}C_L$ \quad
($\mu \notin \frac12 \bZ$; or $\mu \in \frac12+\bZ$ and $\lambda \neq -3,-1,1$; or
$\mu \in \bZ$ and $\lambda \neq -1$);
\item[$(2)$] the cocycle of $(1)$ together with
$\omega(L_m,Y_{n+\frac12})=\delta_{m+n+\mu+\frac12,0}C_{LY}$ \quad
($\mu \in \frac12+\bZ$, $\lambda=-3$);
\item[$(3)$] the cocycle of $(1)$ together with
$\omega(L_m,Y_{n+\frac12})=\frac{m(m-1)}{2}\delta_{m+n+\mu+\frac12,0}C_{LY}$ and
$\omega(M_m,Y_{n+\frac12})=\delta_{m+n+3\mu+\frac12,0}C_{MY}$ \quad
($\mu \in \frac12+\bZ$, $\lambda=-1$);
\item[$(5)$] the cocycle of $(1)$ together with
$\omega(Y_{m+\frac12},Y_{n+\frac12})=-(m+\mu+\frac12)\delta_{m+n+2\mu+1,0}C_Y$ \quad
($\mu \in \bZ$, $\lambda=-1$).
\end{enumerate}
The remaining case $\widetilde{L_{1,\mu}}^{\,4}$ ($\mu \in \frac12+\bZ$, $\lambda=1$) is
discussed in Proposition \ref{prop:l4jac} and Definition \ref{def:l4} below.
\end{definition}

\begin{definition}\label{def:dsvshift}
For $t \in \bZ$, $\alpha,\beta,\gamma,\delta \in \bF$, a linear map $f$ of
$\widetilde{L_{\lambda,\mu}}^{\,i}$ into its centre $Z$ and elements $z_s \in Z$, the
\emph{graded shift of degree $t$} is the linear map $\varphi$ with
$$\varphi(L_m)=\alpha L_{m+t}+\beta M_{m+t}+f(L_m), \qquad
\varphi(M_m)=\gamma M_{m+t}+f(M_m),$$
$$\varphi(Y_{m+\frac12})=\delta Y_{m+t+\frac12}+f(Y_{m+\frac12}), \qquad \varphi(C_s)=z_s .$$
\end{definition}

\begin{theorem}\label{thm:dsv}
Let $i \in \{1,2,3,5\}$ and let $\varphi$ be a graded shift of degree $t$ of
$\widetilde{L_{\lambda,\mu}}^{\,i}$.  Then $\varphi$ is a $\Hom$-Lie structure if and only if
$$\gamma=\delta=\alpha \qquad \text{and} \qquad \alpha\,t=0 .$$
In particular the coefficient $\beta$ is arbitrary, so that
$\mathcal{H}(\widetilde{L_{\lambda,\mu}}^{\,i})$ is infinite-dimensional modulo the maps with
central image.
\end{theorem}

\begin{remark}\label{rem:dsvhalf}
By \cite{Sh24}, $\widetilde{L_{\lambda,\mu}}^{\,i}$ has no non-trivial $\frac12$-derivation for
$i \in \{2,3,5\}$ and for $i=1$ with $\lambda \neq 1$, while
$\Delta(\widetilde{L_{1,\mu}}^{\,1})$ consists of the maps
$$L_m \mapsto \lambda_1 L_m+\sum_t \alpha_t M_{m+t}+\sum_t \beta_t Y_{m+t+\frac12}, \quad
Y_{m+\frac12}\mapsto \lambda_1Y_{m+\frac12}+\sum_t\beta_tM_{m+t+1}, \quad
M_m \mapsto \lambda_1 M_m .$$
The maps with $\beta_t \neq 0$ are homogeneous of half-integer degree and are therefore not
graded shifts in the sense of Definition \ref{def:dsvshift}; by Theorem \ref{thm:half} they are
$\Hom$-Lie structures as well.  The maps with $\alpha_t \neq 0$ are linear combinations of the
graded shifts $\psi_{0,t}\colon L_m\mapsto M_{m+t}$, $M\mapsto 0$, $Y\mapsto 0$ (Definition
\ref{def:dsvshift} with $\alpha=\gamma=\delta=0$, $\beta=1$ and $f=0$), which Theorem
\ref{thm:dsv} exhibits as $\Hom$-Lie structures for \emph{all} $\lambda$ and $\mu$.
\end{remark}

It remains to consider the fifth central extension.  As printed in \cite[Section 4]{Sh24}, the
algebra $\widetilde{L_{1,\mu}}^{\,4}$ ($\mu \in \frac12+\bZ$) is given by
\begin{align}
[L_m,L_n]&=(n-m)L_{m+n}+\tfrac{m^3-m}{12}\delta_{m+n,0}C_L, \label{l4a}\\
[L_m,Y_{n+\frac12}]&=\Big(n+\tfrac12-m+\mu\Big)Y_{m+n+\frac12}
-m(m^2-1)\delta_{m+n+\mu+\frac12,0}C_{LY}, \label{l4b}\\
[L_m,M_n]&=(n-m+2\mu)M_{m+n}-m(m^2-1)\delta_{m+n+2\mu,0}C_M, \label{l4c}\\
[Y_{m+\frac12},Y_{n+\frac12}]&=(n-m)M_{m+n+1}
-(m+\mu)\big((m+\mu)^2-1\big)\delta_{m+n+2\mu,0}C_M . \label{l4d}
\end{align}

\begin{proposition}\label{prop:l4jac}
The multiplication (\ref{l4a})--(\ref{l4d}) does not satisfy the Jacobi identity.  For
$\mu=\frac12$,
$$\big[L_0,[Y_{\frac12},Y_{-\frac12}]\big]+\big[Y_{\frac12},[Y_{-\frac12},L_0]\big]
+\big[Y_{-\frac12},[L_0,Y_{\frac12}]\big]=-\tfrac38\,C_M \neq 0 .$$
More precisely, the two central terms of (\ref{l4c}) and (\ref{l4d}) are supported on pairs of
different total degree: with $\deg L_m=m$, $\deg M_n=n$ and $\deg Y_{n+\frac12}=n+\frac12$, the
first is supported in degree $-2\mu$ and the second in degree $-2\mu+1$.  For such a pair of
supports the cocycle identity (\ref{cocycle}) for the triples $(L,Y,Y)$ forces the central term
of (\ref{l4c}) to vanish identically.
\end{proposition}

The correct central extension is obtained by shifting the support of the $(Y,Y)$-term by one,
and it has a conceptual description.  For $\lambda=1$ the two modules occurring in
$L_{1,\mu}$ are copies of the adjoint module of the Witt part: since $\mu \in \frac12+\bZ$, we
have $2\mu \in \bZ$ and $\mu+\frac12 \in \bZ$, and
$$M_n \longmapsto L_{n+2\mu}, \qquad Y_{n+\frac12}\longmapsto L_{n+\frac12+\mu}$$
transform the brackets $[L_m,M_n]=\big((n+2\mu)-m\big)M_{m+n}$,
$[L_m,Y_{n+\frac12}]=\big((n+\tfrac12+\mu)-m\big)Y_{m+n+\frac12}$ and
$[Y_{m+\frac12},Y_{n+\frac12}]=\big((n+\tfrac12+\mu)-(m+\tfrac12+\mu)\big)M_{m+n+1}$ into the
bracket of the Witt part.  Transporting the Virasoro cocycle along these identifications gives:

\begin{definition}\label{def:l4}
Let $\mu \in \frac12+\bZ$.  The algebra $\widetilde{L_{1,\mu}}^{\,4}$ is the central extension of
$L_{1,\mu}$ by $Z=\langle C_L,C_{LY},C_M\rangle$ given by the cocycle
\begin{align*}
\omega(L_m,L_n)&=\delta_{m+n,0}\tfrac{f(m)}{12}C_L, &
\omega(L_m,Y_{n+\frac12})&=-f(m)\,\delta_{m+n+\mu+\frac12,0}C_{LY},\\
\omega(L_m,M_n)&=-f(m)\,\delta_{m+n+2\mu,0}C_M, &
\omega(Y_{m+\frac12},Y_{n+\frac12})&=-f\big(m+\mu+\tfrac12\big)\,
\delta_{m+n+2\mu+1,0}C_M,
\end{align*}
where $f(x)=x^3-x$.  This is (\ref{l4a})--(\ref{l4d}) with the last line corrected; the
coefficient $-m(m^2-1)=-f(m)$ of (\ref{l4b}) and (\ref{l4c}) is unchanged.
\end{definition}

\begin{proposition}\label{prop:l4lie}
The cocycle of Definition \ref{def:l4} satisfies (\ref{cocycle}), so that
$\widetilde{L_{1,\mu}}^{\,4}$ is a Lie algebra.
\end{proposition}

\begin{theorem}\label{thm:l4}
Let $\mu \in \frac12+\bZ$ and let $\varphi$ be a graded shift of degree $t$ of
$\widetilde{L_{1,\mu}}^{\,4}$ as in Definition \ref{def:dsvshift}.  Then $\varphi$ is a
$\Hom$-Lie structure if and only if
$$\gamma=\delta=\alpha, \qquad \alpha\,t=0, \qquad \beta\,(t+2\mu)=0 .$$
\end{theorem}

\begin{remark}\label{rem:l4}
Theorem \ref{thm:l4} is the only place in this paper where a non-trivial graded shift of
non-zero degree survives in a central extension for one single value of the degree: the map
$L_m \mapsto \beta M_{m-2\mu}$, $M \mapsto 0$, $Y \mapsto 0$ is a $\Hom$-Lie structure of
$\widetilde{L_{1,\mu}}^{\,4}$, and no other map $\psi_{0,t}\colon L_m\mapsto M_{m+t}$,
$M\mapsto 0$, $Y\mapsto 0$ is.  The reason is that the cocycle of Definition \ref{def:l4} pairs
$L_m$ with $M_n$ exactly when $m+n+2\mu=0$: the shift by $-2\mu$ is exactly the one which is
compatible with this pairing.
\end{remark}

\section{Novikov--Witt algebras}\label{sec:nov}

Finally we record what the results of Section \ref{sec:witt} say about the (non-Lie) algebras
of \cite{KKhS25}, whose quasi-derivations are described there.

\begin{definition}[\cite{KKhS25}]\label{def:nov}
For $\xi,\eta \in \bF$ and $\theta \in \bZ\setminus\{0\}$, the \emph{Novikov--Witt algebra}
$\mathrm{W}(\xi,\eta,\theta)$ is the algebra with basis $\{W_i\}_{i \in \bZ}$ and product
$W_n \circ W_m=(\xi+m)W_{n+m}+\eta\,W_{n+m+\theta}$.  For $\varrho \in \bF$, the
\emph{admissible Novikov--Witt algebra} $\mathrm{N}_{\varrho}$ has basis $\{W_i\}_{i \in \bZ}$
and product $W_n \circ W_m=(\varrho+m+2n)W_{n+m}$.
\end{definition}

\begin{proposition}\label{prop:nov}
The commutator algebras of $\mathrm{W}(\xi,\eta,\theta)$ and of $\mathrm{N}_{\varrho}$ are both
isomorphic to the Witt algebra $\WW$ of Definition \ref{def:witt}.  Consequently, for all values
of the parameters, the $\Hom$-Lie structures of these commutator algebras are exactly the maps
$$W_i \longmapsto \sum_{j \in \bZ}\alpha_j W_{i+j} \qquad (i \in \bZ)$$
with finitely many non-zero coefficients $\alpha_j$ independent of $i$.
\end{proposition}

\section{The full matrix algebra and the algebra of upper triangular
matrices}\label{sec:matrix}

The $\frac12$-derivations and the transposed Poisson structures of the Lie algebras $M_n(\bF)$
and $T_n(\bF)$ of all, resp. of all upper triangular, $n \times n$ matrices were determined by
Kaygorodov and Khrypchenko \cite{KKhUT}; the second algebra is the incidence algebra of a chain,
and the general case of a finite poset is treated in \cite{KKhInc}.  In this section we compute
the corresponding spaces of $\Hom$-Lie structures.  They turn out to be considerably larger:
$\Delta(M_n(\bF))$ is two-dimensional and $\Delta(T_n(\bF))$ has dimension $n+2$, whereas
$\mathcal{H}(M_n(\bF))$ has dimension $n^2+1$ and $\mathcal{H}(T_n(\bF))$ has dimension
$n(n+1)+5$.

Throughout the section $\bF$ is a field of characteristic zero, $e_{ij}$ ($1\leq i,j \leq n$) are
the matrix units, $\delta=\sum_{i}e_{ii}$ is the identity matrix, $M_n=M_n(\bF)$ is the Lie
algebra of all $n \times n$ matrices under the commutator, $T_n=T_n(\bF)=\langle e_{ij}\rangle_{i
\leq j}$ is its subalgebra of upper triangular matrices and
$N_n=\langle e_{ij}\rangle_{i<j}$ is the ideal of strictly upper triangular matrices.

\subsection{The full matrix algebra}

\begin{theorem}\label{thm:gl}
Let $\bF$ be an algebraically closed field of characteristic zero, let $n \geq 2$ and let
$Z=\langle \delta \rangle$ be the centre of $M_n$.  Then a linear map $\varphi$ on $M_n$ is a
$\Hom$-Lie structure if and only if
$$\varphi(Z)\subseteq Z \qquad \text{and} \qquad
\pi \circ \varphi|_{\mathfrak{sl}_n} \in \mathcal{H}(\mathfrak{sl}_n),$$
where $\pi \colon M_n \to \mathfrak{sl}_n$ is the projection along $Z$.  Consequently:
\begin{enumerate}
\item if $n \geq 3$, then
$\mathcal{H}(M_n)=\bF\,\id \oplus \Hom(M_n,Z)$ and $\dim \mathcal{H}(M_n)=n^2+1$;
\item if $n=2$, then $\mathcal{H}(M_2)\cong \mathcal{H}(\sltwo)\oplus \Hom(M_2,Z)$ and
$\dim \mathcal{H}(M_2)=6+4=10$.
\end{enumerate}
\end{theorem}

\begin{remark}\label{rem:glhalf}
By \cite[Proposition 13]{KKhUT} the space of $\frac12$-derivations of $M_n$ is
$\Delta(M_n)=\langle \id, \gamma\rangle$, where $\gamma(e_{ij})=\delta_{ij}\delta$, and by
\cite[Theorem 14]{KKhUT} the only non-trivial transposed Poisson structure on $M_n$ is, up to
isomorphism, $e_{ii}\cdot e_{jj}=\delta$.  Both maps $\id$ and $\gamma$ appear in Theorem
\ref{thm:gl}: $\gamma \in \Hom(M_n,Z)$.  The theorem shows how far the inclusion
$\Delta \subseteq \mathcal{H}$ of Theorem \ref{thm:half} is from being an equality here:
$\dim \Delta(M_n)=2$ while $\dim \mathcal{H}(M_n)=n^2+1$.
\end{remark}

\subsection{Upper triangular matrices}

We first record the elementary computation on which everything rests.

\begin{lemma}\label{lem:centmat}
Let $x=\sum_{i \leq j}x_{ij}e_{ij}\in T_n$ and let $1 \leq k<m \leq n$.  Then
$$[x,e_{km}]=\sum_{i \leq k}x_{ik}e_{im}-\sum_{j \geq m}x_{mj}e_{kj},$$
and consequently $[x,e_{km}]=0$ if and only if
$$x_{ik}=0 \ \ (i<k), \qquad x_{mj}=0 \ \ (j>m), \qquad x_{kk}=x_{mm}.$$
\end{lemma}

In particular $\mathcal{C}(T_n)=\langle \delta, e_{1n}\rangle$: a matrix commuting with all of
$N_n$ has, by Lemma \ref{lem:centmat} applied to all $k<m$, all its strictly upper entries equal
to zero except possibly the entry $(1,n)$, and all its diagonal entries equal.

\begin{definition}\label{def:tnmaps}
Let $n \geq 3$.  Define four linear maps on $T_n$ by
$$\begin{array}{llll}
A \colon & e_{12}\mapsto e_{2n}, &
B \colon & e_{11}\mapsto e_{2n}, \ e_{22}\mapsto -e_{2n},\\
A' \colon & e_{n-1,n}\mapsto e_{1,n-1}, \qquad &
B' \colon & e_{nn}\mapsto e_{1,n-1}, \ e_{n-1,n-1}\mapsto -e_{1,n-1},
\end{array}$$
all the remaining basis vectors being sent to $0$.
\end{definition}

\begin{proposition}\label{prop:tnlower}
Let $n \geq 3$.  Then $\id$, $A$, $B$, $A'$, $B'$ and every linear map
$T_n \to \mathcal{C}(T_n)$ are $\Hom$-Lie structures of $T_n$, and the sum
$$\Hom(T_n,\mathcal{C}(T_n))+\langle \id, A,B,A',B'\rangle$$
is direct.  In particular $\dim \mathcal{H}(T_n)\geq n(n+1)+5$.
\end{proposition}

The following lemmas show that these maps exhaust $\mathcal{H}(T_n)$.  The first of them uses
only triples of pairwise commuting elements.

\begin{lemma}\label{lem:tncomm}
Let $\varphi \in \mathcal{H}(T_n)$, let $u=e_{pq}$ with $p \leq q$ and let $k<m$ be such that
$k \neq q$, $m \neq p$ and $l \notin \{p,q\}$ for at least one $l$ with $k \leq l \leq m$.  Then
$$[\varphi(e_{pq}),e_{km}]=0 .$$
\end{lemma}

\begin{lemma}\label{lem:tnstep1}
Let $n \geq 4$, let $\varphi \in \mathcal{H}(T_n)$ and write
$\varphi(e_{pq})=\sum_{i \leq j}x^{pq}_{ij}e_{ij}$.  Then, for all $p \leq q$, the diagonal
entries $x^{pq}_{ii}$ are all equal, with the single exception of $x^{pq}_{pp}$ when $p=q$, of
$x^{12}_{22}$ and of $x^{n-1,n}_{n-1,n-1}$; and the strictly upper entries $x^{pq}_{ij}$
($i<j$) vanish for all $(i,j)$ outside the set $U(p,q)$, where
$$\begin{array}{@{}l@{\ \ }l@{}}
U(p,q)=\{(1,q),(1,n),(p,q),(p,n)\} & \text{if } p<q,\ 
(p,q)\notin\{(1,2),(n-1,n)\},\\
U(1,2)=\{(1,2),(1,n),(2,n)\}, & U(n-1,n)=\{(1,n-1),(1,n),(n-1,n)\},\\
U(p,p)=\{(1,p),(1,n),(p,n)\} \ \ (2\leq p \leq n-1), &\\
U(1,1)=\{(1,n),(2,n)\}, & U(n,n)=\{(1,n-1),(1,n)\} .
\end{array}$$
\end{lemma}

\begin{theorem}\label{thm:tn}
Let $n \geq 4$ and let $\bF$ be a field of characteristic zero.  Then
$$\mathcal{H}(T_n)=\Hom\big(T_n,\mathcal{C}(T_n)\big)\oplus \langle \id, A,B,A',B'\rangle,
\qquad \mathcal{C}(T_n)=\langle \delta, e_{1n}\rangle ,$$
with $A,B,A',B'$ as in Definition \ref{def:tnmaps}.  In particular
$$\dim \mathcal{H}(T_n)=n(n+1)+5 .$$
Explicitly, a linear map $\varphi$ on $T_n$ is a $\Hom$-Lie structure if and only if there are
scalars $s,\beta,\gamma,\beta',\gamma'$ and a linear map $\varrho \colon T_n \to
\langle \delta,e_{1n}\rangle$ such that $\varphi=s\,\id+\varrho$ on all basis vectors except the
six vectors $e_{11},e_{22},e_{12},e_{n-1,n-1},e_{nn},e_{n-1,n}$, on which
$$\begin{array}{ll}
\varphi(e_{11})=s\,e_{11}+\varrho(e_{11})+\gamma e_{2n}, \qquad &
\varphi(e_{nn})=s\,e_{nn}+\varrho(e_{nn})+\gamma' e_{1,n-1},\\
\varphi(e_{22})=s\,e_{22}+\varrho(e_{22})-\gamma e_{2n}, &
\varphi(e_{n-1,n-1})=s\,e_{n-1,n-1}+\varrho(e_{n-1,n-1})-\gamma' e_{1,n-1},\\
\varphi(e_{12})=s\,e_{12}+\varrho(e_{12})+\beta e_{2n}, &
\varphi(e_{n-1,n})=s\,e_{n-1,n}+\varrho(e_{n-1,n})+\beta' e_{1,n-1}.
\end{array}$$
\end{theorem}

\begin{remark}\label{rem:tn3}
For $n=3$ the description of Theorem \ref{thm:tn} remains valid --- note that for $n=3$ the two
exceptional pairs overlap, $e_{22}$ being at the same time $e_{n-1,n-1}$, so that
$\varphi(e_{22})=s\,e_{22}+\varrho(e_{22})-\gamma e_{23}-\gamma' e_{12}$ --- and
$\dim \mathcal{H}(T_3)=3 \cdot 4+5=17$; this case is checked by a direct solution of the linear
system of Corollary \ref{cor:triples} in the $36$ unknown entries of $\varphi$.  For $n=2$ the
algebra $T_2$ is three-dimensional with $\dim T_2^3=1$, so $\dim \mathcal{H}(T_2)=8$ by Theorem
\ref{thm:dim3exact}.
\end{remark}

\begin{remark}\label{rem:tnhalf}
By \cite[Proposition 10]{KKhUT} the $\frac12$-derivations of $T_n$ form the space
$\Delta(T_n)=\langle \id,\alpha\rangle \oplus \langle \beta_i \ : \ 1 \leq i \leq n\rangle$ of
dimension $n+2$, where $\alpha(e_{11})=e_{1n}=-\alpha(e_{nn})$ and $\alpha$ kills all the other
basis vectors, and $\beta_i(e_{jk})=\delta_{(j,k),(i,i)}\delta$.  All of them are visible in
Theorem \ref{thm:tn}: the maps $\beta_i$ lie in $\Hom(T_n,\mathcal{C}(T_n))$, and $\alpha$ lies
in $\Hom(T_n,\mathcal{C}(T_n))$ as well, since $e_{1n}\in \mathcal{C}(T_n)$.  Thus
$\Delta(T_n)\subseteq \Hom(T_n,\mathcal{C}(T_n))+\bF\,\id$, and the four extra maps
$A,B,A',B'$, which are $\Hom$-Lie structures by Proposition \ref{prop:tnlower}, are not
$\frac12$-derivations: for instance $2A([e_{11},e_{12}])=2A(e_{12})=2e_{2n}$, whereas
$[A(e_{11}),e_{12}]+[e_{11},A(e_{12})]=[0,e_{12}]+[e_{11},e_{2n}]=0$.  Again
$\dim \Delta(T_n)=n+2$ is much smaller than $\dim \mathcal{H}(T_n)=n(n+1)+5$.
\end{remark}

\begin{remark}\label{rem:tninc}
The algebra $T_n$ is the incidence algebra $I(X,\bF)$ of the chain
$X=\{1<2<\dots<n\}$, whose $\frac12$-derivations and transposed Poisson structures --- for an
arbitrary finite connected poset $X$ --- were described by Kaygorodov and Khrypchenko
\cite{KKhInc}.  Theorem \ref{thm:tn} is thus the $\Hom$-Lie analogue of \cite{KKhUT,KKhInc} in
the case of a chain.  Lemma \ref{lem:tncomm} holds verbatim for an arbitrary finite connected
poset and yields, for every incidence algebra, the inclusion
$\Hom(I(X,\bF),\mathcal{C}(I(X,\bF)))\subseteq \mathcal{H}(I(X,\bF))$ together with the
constraints of Lemma \ref{lem:centmat}; the determination of $\mathcal{H}(I(X,\bF))$ for a
general poset seems to be an interesting question.
\end{remark}

\section{The Virasoro-like algebra and its $q$-analogue}\label{sec:vlike}

The Virasoro-like algebra is the rank two analogue of the Witt algebra: it is the Lie algebra of
skew derivations of the Laurent polynomial ring in two variables.  Its $\frac12$-derivations, and
more generally its $\delta$-derivations, were recently determined by Lin, Liu and Jiang
\cite{LLJ2}, who proved that non-trivial $\delta$-derivations exist only for $\delta=1$ and that
the algebra admits no non-trivial transposed $\delta$-Poisson structure; the same authors treated
the $q$-analogue in \cite{LLJ1}.  We show that in both cases the $\Hom$-Lie picture is as
rigid as it can be: all $\Hom$-Lie structures are trivial.  This is in sharp contrast with the
Witt algebra (Theorem \ref{thm:witt}), whose $\Hom$-Lie structures form an infinite-dimensional
space, and it gives, by Proposition \ref{prop:homtr}, a new proof of the triviality of the
transposed Poisson structures of these algebras.

For $m=(m_1,m_2)$ and $n=(n_1,n_2)$ in $\bZ^2$ we write
$$\langle m,n \rangle = m_1n_2-m_2n_1 ,$$
so that $\langle \cdot,\cdot\rangle$ is a skew-symmetric bilinear form on $\bZ^2$ with
$\langle m,n\rangle=0$ if and only if $m$ and $n$ are proportional over $\mathbb{Q}$.

\begin{definition}[\cite{LLJ2}]\label{def:vlike}
The \emph{Virasoro-like algebra} $\mathcal{V}$ is the Lie algebra with basis
$\{L_m\}_{m \in \bZ^2\setminus \{0\}}$ and multiplication
$$[L_m,L_n]=\langle m,n\rangle L_{m+n}, \qquad m,n \in \bZ^2\setminus\{0\},$$
with the convention $L_0=0$.
\end{definition}

The convention $L_0=0$ is consistent: if $m+n=0$, then $\langle m,n\rangle=0$ anyway.  The
algebra of \cite{LLJ2}, where the structure constant is $\langle n,m\rangle$, is obtained from
$\mathcal{V}$ through the isomorphism $L_m \mapsto -L_m$.

\begin{theorem}\label{thm:vlike}
All the $\Hom$-Lie structures of the Virasoro-like algebra are trivial:
$$\mathcal{H}(\mathcal{V})=\bF\,\id .$$
\end{theorem}

\begin{corollary}\label{cor:vliketp}
The Virasoro-like algebra admits no non-trivial transposed Poisson structure.
\end{corollary}

\begin{remark}\label{rem:vlikellj}
Corollary \ref{cor:vliketp} is the case $\delta=\frac12$ of \cite[Theorem 3.5]{LLJ2}; here it
is deduced from a stronger statement, namely that the much larger space
$\mathcal{H}(\mathcal{V})\supseteq \Delta(\mathcal{V})$ is already reduced to the scalars.  It is
worth stressing how different the rank one and the rank two situations are.  For the Witt algebra
$\WW$ every ``shift with constant coefficients'' $e_i \mapsto e_{i+t}$ is a $\Hom$-Lie structure
(Theorem \ref{thm:witt}), whereas for $\mathcal{V}$ the analogous map $L_m \mapsto L_{m+t}$
satisfies, by a direct computation,
$$E_{\varphi}(L_a,L_b,L_c)=-\big(\langle a,b\rangle+\langle b,c\rangle+\langle c,a\rangle\big)
\langle a+b+c,t\rangle\,L_{a+b+c+t},$$
which is non-zero as soon as $t \neq 0$.  The two algebras also behave differently with respect
to $\frac12$-derivations.  Every shift $e_i \mapsto e_{i+t}$ of $\WW$ is a
$\frac12$-derivation, since
$[e_{i+t},e_j]+[e_i,e_{j+t}]=\big((i+t-j)+(i-j-t)\big)e_{i+j+t}=2(i-j)e_{i+j+t}$, so that by
Theorem \ref{thm:witt} and Theorem \ref{thm:half} one has $\Delta(\WW)=\mathcal{H}(\WW)$, both
being the infinite-dimensional space of shifts with constant coefficients.  For the
Virasoro-like algebra, on the contrary, $\Delta(\mathcal{V})=\mathcal{H}(\mathcal{V})=\bF\,\id$
by Theorem \ref{thm:vlike} and \cite{LLJ2}.
\end{remark}

We now turn to the $q$-analogue.  Let $q \in \bF^{\times}$.

\begin{definition}[\cite{LLJ1}]\label{def:qvlike}
The \emph{$q$-analogue of the Virasoro-like algebra} is the Lie algebra $\mathcal{V}_q$ with
basis $\{L_m\}_{m \in \bZ^2\setminus\{0\}}$, the convention $L_0=0$, and multiplication
$$[L_m,L_n]=\lambda(m,n)L_{m+n}, \qquad \lambda(m,n)=q^{m_2n_1}-q^{m_1n_2} .$$
\end{definition}

If $q$ is not a root of unity, then
\begin{equation}\label{lambdazero}
\lambda(m,n)=0 \iff m_2n_1=m_1n_2 \iff \langle m,n\rangle =0 ,
\end{equation}
i.e. $\lambda$ vanishes exactly where the form $\langle \cdot,\cdot\rangle$ of the classical case
does.  Theorem \ref{thm:qvlike} below rests on this fact, together with the following identity.

\begin{lemma}\label{lem:qlambda}
Let $a,b \in \bZ^2$ and put $u=q^{a_2b_1}$, $v=q^{a_1b_2}$, $\alpha=q^{a_1a_2}$,
$\beta=q^{b_1b_2}$.  Then
$$\lambda(b,a+b)=\beta(v-u), \ \ \lambda(a+b,a)=\alpha(v-u), \ \
\lambda(a,a+2b)=\alpha(u^2-v^2), \ \ \lambda(b,2a+b)=\beta(v^2-u^2).$$
In particular $\lambda(b,a+b)\lambda(a,a+2b)=-\lambda(a+b,a)\lambda(b,2a+b)
=-\alpha\beta(v-u)^2(v+u)$.
\end{lemma}

\begin{theorem}\label{thm:qvlike}
Let $q \in \bF^{\times}$ not be a root of unity.  Then all the $\Hom$-Lie structures of
$\mathcal{V}_q$ are trivial:
$$\mathcal{H}(\mathcal{V}_q)=\bF\,\id .$$
Consequently $\mathcal{V}_q$ admits no non-trivial transposed Poisson structure.
\end{theorem}

\begin{remark}\label{rem:qvlike}
By \cite{LLJ1} the algebra $\mathcal{V}_q$ with $q$ generic has no non-trivial
$\frac12$-derivation, and hence no non-trivial transposed Poisson structure; Theorem
\ref{thm:qvlike} strengthens this to the statement that $\mathcal{V}_q$ has no non-trivial
$\Hom$-Lie structure at all.  For $q$ a primitive root of unity the situation changes:
\cite{LLJ1} exhibits non-trivial $\frac12$-derivations of $\mathcal{V}_q$, which are then
non-trivial $\Hom$-Lie structures by Theorem \ref{thm:half}, so that
$\mathcal{H}(\mathcal{V}_q)\neq \bF\,\id$ in that case.  The hypothesis that $q$ is not a
root of unity enters only through the equivalence (\ref{lambdazero}) and the non-vanishing of
$u+v$ in Lemma \ref{lem:qlambda}.
\end{remark}

\section{The planar Galilean conformal algebra}\label{sec:pgca}

Wu and Zhang \cite{WZ} proved that every $\frac12$-derivation of the planar Galilean conformal
algebra is a scalar, so that all its transposed Poisson structures are trivial.  We show that the
algebra nevertheless carries an infinite-dimensional space of $\Hom$-Lie structures, all of
which, apart from the scalars, map the Witt part into the abelian part.

\begin{definition}[\cite{WZ}]\label{def:pgca}
Let $\Gamma$ be a non-trivial additive subgroup of $\bF$.  The \emph{planar Galilean conformal
algebra} $\mathcal{G}(\Gamma)$ is the Lie algebra with basis
$\{L_{\mu},H_{\mu},I_{\mu},J_{\mu}\}_{\mu \in \Gamma}$ and multiplication
$$\begin{array}{lll}
{}[L_{\mu},L_{\nu}]=(\mu-\nu)L_{\mu+\nu}, \quad &
{}[L_{\mu},H_{\nu}]=-\nu H_{\mu+\nu}, \quad &
{}[L_{\mu},I_{\nu}]=(\mu-\nu)I_{\mu+\nu},\\
{}[L_{\mu},J_{\nu}]=(\mu-\nu)J_{\mu+\nu}, &
{}[H_{\mu},I_{\nu}]=J_{\mu+\nu}, &
{}[H_{\mu},J_{\nu}]=-I_{\mu+\nu},
\end{array}$$
all the remaining brackets of basis vectors being zero.  For $\Gamma=\bZ$ this is the algebra of
\cite{WZ}.
\end{definition}

The subspace $\mathcal{I}=\langle I_{\mu},J_{\mu}\rangle_{\mu \in \Gamma}$ is an abelian ideal,
and $\langle L_{\mu}\rangle_{\mu}$ is a copy of the Witt algebra (of the higher rank Witt algebra
$W(\Gamma)$ in general); the subalgebra $\langle L_{\mu},H_{\mu}\rangle_{\mu}$ is the
Heisenberg--Virasoro algebra $W(\Gamma;0,0)$ of Definition \ref{def:wab}.  Since
$\mathcal{G}(\Gamma)$ is $\Gamma$-graded, with
$\mathcal{G}(\Gamma)_{\mu}=\langle L_{\mu},H_{\mu},I_{\mu},J_{\mu}\rangle$, we look, as in
Sections \ref{sec:wab}--\ref{sec:nfg}, at the homogeneous maps.

\begin{definition}\label{def:pgcashift}
For $t \in \Gamma$ and a matrix $(c_{XY})_{X,Y \in \{L,H,I,J\}}$ of scalars, the \emph{graded
shift} of degree $t$ with coefficients $c_{XY}$ is the linear map $\varphi$ with
$$\varphi(X_{\mu})=\sum_{Y \in \{L,H,I,J\}} c_{XY}\,Y_{\mu+t} \qquad (\mu \in \Gamma).$$
\end{definition}

\begin{theorem}\label{thm:pgca}
A graded shift $\varphi$ of degree $t$ of $\mathcal{G}(\Gamma)$ is a $\Hom$-Lie structure if and
only if
$$c_{LL}=c_{HH}=c_{II}=c_{JJ}=:\alpha, \qquad \alpha\,t=0,$$
and all the remaining coefficients vanish except $c_{LI}$ and $c_{LJ}$, which are arbitrary.  In
other words, the graded $\Hom$-Lie structures of $\mathcal{G}(\Gamma)$ are exactly the maps
$$L_{\mu}\mapsto \alpha L_{\mu+t}+\beta I_{\mu+t}+\gamma J_{\mu+t}, \qquad
H_{\mu}\mapsto \alpha H_{\mu+t}, \qquad I_{\mu}\mapsto \alpha I_{\mu+t}, \qquad
J_{\mu}\mapsto \alpha J_{\mu+t},$$
with $\beta,\gamma \in \bF$ arbitrary and $\alpha=0$ unless $t=0$.
\end{theorem}

\begin{corollary}\label{cor:pgca}
For all $t \in \Gamma$ and $\beta,\gamma \in \bF$ the map
$$\psi_{t,\beta,\gamma} \colon \quad L_{\mu}\mapsto \beta I_{\mu+t}+\gamma J_{\mu+t}, \qquad
H_{\mu},I_{\mu},J_{\mu}\mapsto 0$$
is a $\Hom$-Lie structure of $\mathcal{G}(\Gamma)$, and so is every finite sum
$\alpha\,\id+\sum_{t \in S}\psi_{t,\beta_t,\gamma_t}$.  Hence
$\mathcal{H}(\mathcal{G}(\Gamma))$ is infinite-dimensional whenever $\Gamma$ is infinite, while
by \cite[Theorem 3.7]{WZ} all the $\frac12$-derivations and all the transposed Poisson
structures of $\mathcal{G}(\bZ)$ are trivial.
\end{corollary}

\begin{remark}\label{rem:pgcawitt}
Theorem \ref{thm:pgca} says in particular that no graded shift of non-zero degree acting by a
non-zero scalar on the $L$'s is a $\Hom$-Lie structure of $\mathcal{G}(\Gamma)$: the Witt-type
shifts $L_{\mu}\mapsto L_{\mu+t}$, which are $\Hom$-Lie structures of the Witt algebra by
Theorem \ref{thm:witt}, do not survive the extension by $H$, $I$ and $J$.  This is the same
phenomenon as in Remark \ref{rem:wabhv} for the Heisenberg--Virasoro algebra $W(\Gamma;0,0)$,
which is a subalgebra of $\mathcal{G}(\Gamma)$; in fact the obstruction already appears on the
triples $(L_{\mu},L_{\nu},H_{\xi})$, whose $H$-component is
$(\mu-\nu)\big(c_{HH}(\xi+t)-c_{LL}\xi\big)$ both here and there.
\end{remark}

\section{Further directions}\label{sec:further}

The methods used above apply to many further classes of Lie algebras, which seem to deserve a
systematic study.

\begin{enumerate}
\item The complete classification --- beyond the graded maps described in Theorem \ref{thm:wab}
--- of the $\Hom$-Lie structures of the algebras $W(\Gamma;a,b)$, and in particular of the
Heisenberg--Virasoro algebra $W(0,0)$.
\item The complete classification --- beyond the graded shifts --- of the $\Hom$-Lie structures
of the central extensions of Sections \ref{sec:cext}--\ref{sec:dsv}, and the determination of
the $\Hom$-Lie structures of the central extensions of $HW_{-1}(G)$ (see Proposition
\ref{prop:hwjac}, where the table printed in \cite{KKhS24b} is shown not to be a Lie algebra).
\item The Galilean conformal algebras, which contain $\WW$ together with a non-trivial module;
their transposed Poisson structures are described in \cite{KLZ23}.  For the planar Galilean
conformal algebra the graded $\Hom$-Lie structures are determined in Theorem \ref{thm:pgca};
the description of the non-graded ones remains open.
\item The Block algebras, the algebras of Weyl type, and the higher rank Witt algebras
$W(\Gamma)$ for an additive subgroup $\Gamma$ of the ground field, whose transposed Poisson
structures are known by \cite{KKh23,KL23}.
\item Loop (current) algebras $\mathfrak{g} \otimes \bF[t,t^{-1}]$ of a simple finite-dimensional
algebra $\mathfrak{g}$, and affine Kac--Moody algebras.
\item Direct, elementary proofs --- in the spirit of Corollary \ref{cor:triples} and Theorem
\ref{thm:sl2}, that is, by an analysis of the triples of root vectors --- of the triviality of
the $\Hom$-Lie structures of the finite-dimensional simple Lie algebras
$\mathfrak{g}=\mathfrak{sl}_n$, $\mathfrak{so}_n$, $\mathfrak{sp}_{2n}$ with
$\mathfrak{g}\not\cong \sltwo$, which is known by
\cite[Theorem 3.3]{xjl15}, and of their analogues in prime characteristic.
\item Nilpotent Lie algebras of nilpotency class at least $3$, filiform algebras for instance,
where Proposition \ref{prop:2step} no longer applies; and solvable Lie algebras with non-abelian
nilpotent radical, extending Theorem \ref{thm:codim1}; for the transposed Poisson counterpart of
this question see \cite{KKhu24}.
\item The $\Hom$-Lie structures of the incidence algebra $I(X,\bF)$ of an arbitrary finite
connected poset $X$, of which Theorem \ref{thm:tn} treats the case of a chain; and those of the
$q$-analogue $\mathcal{V}_q$ of the Virasoro-like algebra when $q$ is a primitive root of unity,
a case in which non-trivial $\frac12$-derivations do exist \cite{LLJ1} and Theorem
\ref{thm:qvlike} does not apply.
\item The Schr\"odinger algebra $\mathcal{S}_n$ in $(n+1)$-dimensional space-time and the
quasi-filiform Lie algebras of maximum length, whose $\frac12$-derivations are described in
\cite{YTK,ADSS}; both are finite-dimensional, so that the bounds of Corollary \ref{cor:dimn} and
the method of Lemma \ref{lem:tncomm} apply to them.
\item The determination, for each $n$, of the minimum of $\dim \mathcal{H}(\LL)$ over all
$n$-dimensional Lie algebras $\LL$, and of the classes of algebras for which the bounds of
Corollary \ref{cor:dimn} are attained.  By Theorems \ref{thm:dim3exact} and \ref{thm:sl2} this
minimum equals $6$ for $n=3$, attained exactly on $\sltwo$;
Examples \ref{ex:fil4} and \ref{ex:sl2ab} show that the refined bound of Theorem \ref{thm:dimn}
is attained in dimension four, and that the bound of Proposition \ref{prop:centralizer} is
asymptotically sharp.
\end{enumerate}

\section*{Acknowledgements}

The authors are grateful to Ivan Kaygorodov for the idea of this paper and for his advice
during its writing. We also thank Kseniia Viatkina for her work on an early version
of this paper.

\bigskip

\noindent \textbf{Jobir Adashev}\\
\textit{Affiliation}: Institute of Mathematics, Academy of Sciences of Uzbekistan, Tashkent,
Uzbekistan\\
\textit{Email address}: adashevjq@mail.ru

\medskip

\noindent \textbf{Majidkhon Azizov}\\
\textit{Affiliation}: Institute of Mathematics, Academy of Sciences of Uzbekistan, Tashkent,
Uzbekistan\\
\textit{Email address}: azizovmajidkhan@gmail.com

\medskip

\noindent \textbf{Vasily Voronin}\\
\textit{Affiliation}: Sobolev Institute of Mathematics, Novosibirsk, Russia\\
\textit{Email address}: voronin.vasily@gmail.com

\end{document}